\documentclass[a4paper]{amsart}
\usepackage[utf8]{inputenc}
\usepackage[english]{babel}
\usepackage{float}

\usepackage[normalem]{ulem}
\usepackage{amscd,amssymb}
\usepackage{amsmath}
\usepackage{amsthm}
\usepackage{graphicx}
\usepackage{fancybox}
\usepackage{epic,eepic}
\usepackage{amstext}
\usepackage{mathtools}
\usepackage{esint}
\usepackage[square,numbers]{natbib}
\usepackage{booktabs}
\usepackage{multirow}
\usepackage[hidelinks]{hyperref}
\newcommand{\weblink}[2]{%
  \href{#1}{\textsf{#2}\,\textsuperscript{\tiny$\nearrow$}}%
}
\usepackage{pgfplots}
\pgfplotsset{compat=1.18}
\usepackage{subcaption}
\usepackage[noabbrev, capitalise, nameinlink]{cleveref}
\usepackage{aliascnt}

\newtheorem{theorem}{Theorem}[section]

\newaliascnt{proposition}{theorem}

\aliascntresetthe{proposition}
\crefname{proposition}{proposition}{propositions}
\Crefname{proposition}{Proposition}{Propositions}

\newaliascnt{example}{theorem}
\newtheorem{example}[example]{Example}
\aliascntresetthe{example}
\crefname{example}{example}{examples}
\Crefname{example}{Example}{Examples}

\theoremstyle{definition}
\newaliascnt{definition}{theorem}
\newtheorem{definition}[definition]{Definition}
\aliascntresetthe{definition}

\crefname{definition}{definition}{definitions}
\Crefname{definition}{Definition}{Definitions}

\theoremstyle{remark}
\newtheorem{remark}[theorem]{Remark}
\numberwithin{theorem}{section}
\numberwithin{equation}{section}
\numberwithin{figure}{section}

\newcommand{\calD}{\mathcal{D}}
\newcommand{\calE}{\mathcal{E}}
\newcommand{\calF}{\mathcal{F}}

\newcommand{\calH}{\mathcal{H}}
\newcommand{\calM}{\mathcal{M}}

\newcommand{\frakP}{\mathfrak{P}}
\newcommand{\frakQ}{\mathfrak{Q}}

\newcommand{\calY}{\mathcal{Y}}
\newcommand{\calZ}{\mathcal{Z}}
\newcommand{\R}{\mathbb{R}}
\newcommand{\N}{\mathbb{N}}
\newcommand{\dd}{\textrm{d}}

\newcommand{\la}{\langle}
\newcommand{\ra}{\rangle}

\DeclareMathOperator{\Ima}{Im}

\DeclareMathOperator*{\ddiv}{div}
\DeclareMathOperator*{\spn}{span}
\DeclareMathOperator*{\argmin}{arg\,min}

\newcommand{\wLapA}{\mathfrak{L}_A}

\newcommand{\pmDis}{\varrho}
\newcommand{\pmFun}{v}
\newcommand{\pmComp}{\Gamma} %
\newcommand{\resDis}{b}
\newcommand{\resFun}{r}

\newcommand{\jacDis}{J}

\newcommand{\paramLinear}{c}
\newcommand{\paramPM}{\theta}
\newcommand{\paramUpdate}{\xi}
\newcommand{\linform}{\lambda}
\newcommand{\rieszform}{\zeta}
\newcommand{\testform}{\mathfrak{z}}
\newcommand{\greenlift}{\mathfrak{g}}
\newcommand{\Id}[1][d]{I_#1}

\newcommand{\archChart}{\psi}
\newcommand{\activation}{\chi}
\newcommand{\braket}[2]{\langle #1\,,\, #2\rangle}
\newcommand{\norm}[1][\cdot]{\left\|#1\right\|}

\allowdisplaybreaks

\title[Beyond PINNs: A Unified %
Framework for Neural and Hybrid PDE Solvers]{Beyond PINNs: A Unified Gauss--Newton and Petrov--Galerkin Framework for Neural and Hybrid PDE Solvers}
\author[N.~Schwencke, R.~Maier]{Nilo Schwencke${}^*$ and Roland Maier${}^\dagger$}
\address{${}^*$
ENS de Lyon, CNRS, Université Claude Bernard Lyon 1, Inria, LIP, UMR 5668, 69342, Lyon cedex 07, France}
\email{nilo.schwencke@ens-lyon.fr}
\address{${}^\dagger$
Institute for Applied and Numerical Mathematics, Karlsruhe Institute of Technology, Englerstr.~2, 76131 Karlsruhe, Germany}
\email{roland.maier@kit.edu}

\begin{document}

\begin{abstract}
Physics-informed neural networks and finite element methods provide two
different paradigms for the numerical approximation of partial differential
equations: the former are commonly trained by minimizing pointwise strong
residuals, whereas the latter are naturally built from weak variational
formulations and the finite-dimensional systems obtained after discretization.
In this work, we introduce a common framework based on the discretization of functional Gauss--Newton problems by finite families of linear measurements. We show that, through an appropriate duality pairing, the linear measurements can be represented by test functions. The resulting Gauss--Newton system is then precisely a Petrov--Galerkin discretization of the linearized functional problem. This perspective recovers pointwise collocation and natural-gradient constructions as particular cases, while making the choice of test functions an explicit algorithmic design choice. 
We specialize this framework to elliptic problems, where it naturally leads to
weak residual formulations and to a hybrid finite element--neural construction acting
on complementary approximation spaces.
Numerical experiments support the proposed
framework and demonstrate the effectiveness of weak Gauss--Newton formulations
and hybrid finite element--neural approximations.

\end{abstract}



\maketitle

\section{Introduction}
\label{sec:introduction}

Partial differential equations (PDEs) constitute one of the principal
mathematical tools for modeling physical phenomena across science and
engineering. Their numerical approximation has consequently motivated a wide
range of methods, among which Galerkin and Petrov--Galerkin discretizations
play a central role: the solution is approximated in a prescribed
finite-dimensional approximation space, while the governing equation is
enforced against a suitable family of test functions through its variational
formulation. Finite element methods are one of their most successful
realizations, typically using piecewise-polynomial approximation spaces
constructed on a mesh; see, e.g.,~\cite{Cia78,Bra07,BreS08}.

Neural-network-based PDE solvers provide a rather different approximation
paradigm. In particular, physics-informed neural networks
(PINNs)~\cite{dissanayakeNeuralnetworkbasedApproximationsSolving1994,
lagarisArtificialNeuralNetworks1998,raissiPhysicsinformedNeuralNetworks2019}
represent the solution by a nonlinear parametric model and determine its
parameters by minimizing residuals of the governing equation and boundary
conditions. Their mesh-free character and the flexibility of the neural
approximation model make them attractive in settings where classical
discretizations may be difficult to construct. Although early PINN
formulations often exhibited limited accuracy and difficult optimization,
substantial progress has been obtained through adaptive sampling
\cite{nabianEfficientTrainingPhysicsinformed2021,
wuComprehensiveStudyNonadaptive2023,
maoPhysicsinformedNeuralNetworks2023,
dawMitigatingPropagationFailures2023,
nguyenFixedBudgetOnlineAdaptive2023,
lauPINNACLEPINNAdaptive2024a}
and, in particular, through Gauss--Newton and natural-gradient-based
optimization
\cite{mullerAchievingHighAccuracy2023,SchF25,
jniniGaussnewtonNaturalGradient2025,
mckayNearoptimalSketchyNatural2025,
schwencke2025amstramgram,
jniniDualNaturalGradient2025,
jniniCurvatureAwareOptimizationHighAccuracy2026,
nouyNaturalGradientDescent2026,
mckayErrorWhiteningWhy2026,webbOptimisationFrameworkWellConditioned2026}.

Despite these developments, the classical finite element and PINN viewpoints
remain rather different. Standard PINNs usually start from a strong
formulation and construct a finite residual vector by evaluating the
differential equation at collocation points. For a second-order elliptic
operator, this requires the strong residual and its parameter derivatives to
possess sufficient pointwise regularity. In contrast, the natural formulation
used by finite element methods is typically weak: the elliptic residual belongs
to $H^{-1}(\Omega)$ and is evaluated through its action on test functions in
$H_0^1(\Omega)$. This distinction becomes particularly relevant for
nonsmooth solutions, distributional right-hand sides, or FE
functions themselves, for which a pointwise strong residual may be unnatural
or even undefined.

Several neural approaches have investigated variational or weak formulations.
VPINNs and hp-VPINNs
\cite{kharazmiVariationalPhysicsInformedNeural2019,
kharazmiHpVPINNsVariationalPhysicsinformed2021}
construct variational losses from prescribed test spaces, while related
Petrov--Galerkin methods allow more general finite element or neural test
families \cite{shangDeepPetrovGalerkinMethod2022}.
Other approaches adapt the test spaces themselves, through adversarial optimization~\cite{zangWeakAdversarialNetworks2020} or minimum-residual formulations~\cite{rojasRobustVariationalPhysicsInformed2024} connected to the classical
theory of optimal Petrov--Galerkin test spaces~\cite{demkowiczClassDiscontinuousPetrov2011}.

Closer to the present finite-measurement viewpoint, a related construction has
recently been proposed in \cite{bonStableNonlinearDynamical2025} in the
context of nonlinear dynamical approximation for time-dependent PDEs.
Therein, the parameters of a nonlinear decoder are evolved by projecting the
PDE dynamics onto its tangent space, with this projection approximated through
finite families of linear observations. The analysis focuses in particular
on the stability of the resulting reconstruction and on the adaptive selection
of the observations.
Our setting differs in that the finite measurements arise from the
discretization of a functional Gauss--Newton problem, a perspective that will
lead to the Petrov--Galerkin interpretation developed below.

At the same time, hybrid methods combining finite element and neural
approximations have attracted increasing attention, see, e.g.,
\cite{FraMicNavVig2025,BarDupFauFraLecLleMicVic2025,
margenbergDNNMGHybridNeural2024}. These developments raise a common
question: can the finite residual systems used by modern Gauss--Newton neural
solvers, weak Petrov--Galerkin formulations, and hybrid finite element--neural approximations be cast within a single construction and
handled by the same optimization machinery?
\par\addvspace{0.5\baselineskip}

The starting point of this work is the observation that Gauss--Newton,
formulated at the functional level, provides such a common framework. By
applying a finite family of linear measurements to the linearized functional
problem, one obtains an ordinary finite residual vector and Jacobian to which
standard Gauss--Newton algorithms can be applied.
When these measurements are represented through pairings with test functions,
the resulting system is
precisely a Petrov--Galerkin discretization of the linearized problem. This
interpretation recovers natural-gradient and pointwise collocation
constructions as particular cases, while making the choice of measurements,
and hence of test functions, an independent component of the numerical
method.

For elliptic PDEs, this flexibility naturally leads to weak residuals in
$H^{-1}(\Omega)$ tested against functions in $H_0^1(\Omega)$, allowing the
same Gauss--Newton machinery to operate directly at the variational level with
a broad choice of test families, including but not limited to
operator-adapted Green sections. The same approximation--test viewpoint also
suggests a hybrid finite element--neural construction in which a FE
component captures a prescribed approximation space, while the neural
component and the test functions act on its energy-orthogonal complement.
This yields a unified weak Gauss--Newton formulation without requiring alternating optimization of the finite element and neural components.
\par\addvspace{0.5\baselineskip}

The main contributions of this work are therefore threefold:
\begin{enumerate}
 \item we establish a Petrov--Galerkin interpretation of functional Gauss--Newton problems discretized by
finite linear measurements;
 \item we extend this construction to weak elliptic
residuals with general test families;
 \item we derive a hybrid finite element--neural formulation in which the finite element and neural components act on complementary energy subspaces.
\end{enumerate}
Numerical experiments
validate these constructions on both smooth and low-regularity problems.
They show that standard Gauss--Newton-based solvers can be applied directly to
weak residuals, reaching accuracies comparable to or substantially higher than
specialized energy-based methods,
and that the hybrid finite element--neural formulation can substantially
improve the robustness of generic weak test spaces while remaining competitive with
standalone finite element discretizations with a comparable number of degrees of freedom.

The remainder of the paper is organized as follows.
\Cref{sec:problem-methods-intro} introduces the model PDE, parametric
approximation models, Galerkin and Petrov--Galerkin methods, neural networks,
and PINNs.
\Cref{sec:anagram-green} develops finite-dimensional and functional
Gauss--Newton frameworks, introduces discretization by linear measurements, and derives
its Petrov--Galerkin interpretation together with the weak-test formulation.
\Cref{sec:hybrid-method} presents the hybrid finite element--neural model and its
associated weak-test Gauss--Newton update.
\Cref{sec:numerical-experiments} contains the numerical validation.
Additional connections with natural gradient, kernel approximation,
pointwise collocation, sketching, and regularization are presented in \Cref{app:gn-approximation-and-test-examples}. We also refer to a companion blog (see \url{https://nilo.schwencke.me/tutorials/beyond-pinns-companion/}) for more detailed
experimental protocols and diagnostics as well as additional numerical examples.

\section{Approximation of PDE solutions}
\label{sec:problem-methods-intro}

\subsection{Model problem}
\label{subsec:model-problem}

We consider a Hilbert space $\calH \subset H^1(\Omega)$,
where $\Omega \subset \R^d$ ($d \in \N$) is a bounded Lipschitz domain.\footnote{More generally, one could consider an abstract Hilbert space $\calH$ of functions from $\Omega$ to $\R$, but the present setting is sufficient for the purposes of this article.} We define operators
\begin{equation}\label{eq:ops}
    D\colon \calH \to L^2(\Omega) \qquad\text{and}\qquad B\colon \calH \to L^2(\partial \Omega),
\end{equation}
where $D$ is a linear differential operator and $B$ a linear boundary condition. We define a scalar product $\la \bullet, \bullet \ra \colon \calH \times \calH \to \R$ by
\begin{equation}\label{eq:defScaPro}
    \la v, w \ra \coloneqq (D[v], D[w])_{L^2(\Omega)} + (B[v],B[w])_{L^2(\partial \Omega)},
\end{equation}
possibly factoring out the nullspace of the pair $(D,B)$.
Given $f\in L^2(\Omega)$ and $g_{\mathrm{bc}}\in L^2(\partial\Omega)$,
an abstract PDE formulation now seeks a function $u\in\calH$ such that
\begin{equation}\label{eq:PDEstrong}
\begin{aligned}
    D[u] &= f \;&&\text{in }L^2(\Omega),\\
    B[u] &= g_{\mathrm{bc}} \;&&\text{in }L^2(\partial \Omega).
\end{aligned}
\end{equation}
We refer to~\eqref{eq:PDEstrong} as a \emph{boundary value problem}.

\begin{example}[second-order elliptic PDE]\label{ex:pde}
    A classical example is the boundary value problem
    \begin{equation}\label{eq:ellPDE}
    \begin{aligned}
        -\ddiv(A\nabla u) &= f \; &&\mathrm{in}\;\Omega,\\
        u &= g_{\mathrm{bc}} &&\mathrm{on}\;\partial\Omega,
    \end{aligned}
    \end{equation}
    where $A\in L^\infty(\Omega,\R^{d\times d})$ is symmetric and, for constants $0<\alpha\le\beta$, satisfies
    \begin{equation}\label{eq:propA}
        \alpha|\eta|^2
        \le
        (A(x)\eta)\cdot\eta
        \le
        \beta|\eta|^2
    \end{equation}
    for almost all $x\in\Omega$ and every $\eta\in\R^d$. Here,
    $|\bullet|$ denotes the Euclidean norm in $\R^d$.
    In the strong formulation~\eqref{eq:PDEstrong}, one takes
    $D=-\ddiv(A\nabla\bullet)$ and $B=\gamma$ on a sufficiently regular
    space $\calH_{\mathrm{str}}\subset H^1(\Omega)$ such that
    $D[u]\in L^2(\Omega)$ and the trace $\gamma u$ is well defined.
    For homogeneous Dirichlet data, the corresponding weak formulation is
    naturally posed on $H_0^1(\Omega)$, where the elliptic operator is instead
    understood as a map into $H^{-1}(\Omega) =
    (H_0^1(\Omega))'$.
\end{example}

\subsection{Parametric models}
\label{subsec:parametric-models}

Many approximation methods can be formulated in terms of a finite number of degrees of freedom. This motivates the introduction of an abstract notion that encompasses both classical discretization techniques and modern machine learning models.

Let $p \in \N$. A \emph{parametric model} is a differentiable map
\begin{equation}\label{eq:param}
    v :
    \R^p
    \to
    \calH,
    \qquad
    \paramPM
    \mapsto
    v_\paramPM .
\end{equation}
The vector $\paramPM \in \R^p$ collects the degrees of freedom of the approximation.
As the parameters vary, the model generates a family of functions
\begin{equation}
    \calM_\pmFun
    \coloneqq
    \{v_\paramPM \mid \paramPM\in\R^p\}
    \subset
    \calH,
\end{equation}
which may be viewed, under mild assumptions, as a finite-dimensional manifold embedded in $\calH$. Understanding how the model can evolve therefore amounts to studying the local geometry of $\calM_\pmFun$.
The differential of the parametrization at $\paramPM$ is then given by
\begin{equation}\label{eq:der}
    \dd v_\paramPM :
    \R^p
    \to
    \calH,
    \qquad
    \paramUpdate
    \mapsto
    \sum_{i=1}^{p}
    \paramUpdate_i
    \partial_i v_\paramPM .
\end{equation}
Its image defines the tangent space of the manifold at $v_\paramPM$, given by
\begin{equation}\label{eq:imSpaces}
    T_\paramPM\calM_\pmFun
    \coloneqq
    \Ima(\dd v_\paramPM)
    =
    \spn\{
        \partial_i v_\paramPM
        \mid
        i=1,\ldots,p
    \}.
\end{equation}
The tangent space describes the admissible infinitesimal variations of the
model around $v_\paramPM$.

Before considering particular examples, it is useful to introduce one further
construction. Let $v:\R^p\to\calH$ be a parametric model and let
$\calF:\calH\to\calY$ be a differentiable map between functional spaces. We
call the composition
\begin{equation}
    \pmComp^{\calF}
    \coloneqq
    \calF\circ v:
    \R^p
    \to
    \calY,
    \qquad
    \pmComp_\paramPM^{\calF}
    =
    \calF(v_\paramPM),
    \label{eq:compound-parametric-model}
\end{equation}
a \emph{compound parametric model}.
Note that $\pmComp^{\calF}$ is itself a parametric model, now taking values in $\calY$.
Its admissible first-order variations
are described by
\begin{equation*}
    T_\paramPM\calM_{\pmComp^{\calF}}
    =
    \Ima\bigl(\dd\pmComp_\paramPM^{\calF}\bigr)
    =
    \spn\{
        \partial_i\pmComp_\paramPM^{\calF}
        \mid 1\leq i\leq p
    \}
    \subset\calY.
\end{equation*}

Two important examples of parametric models are finite element approximations
and neural networks, which we introduce next. Compound models will in
particular arise naturally in the PINN formulation below.

\subsection{Galerkin and Petrov--Galerkin methods}
\label{subsec:galerkin-petrov}

To simplify the presentation, within this subsection we stick to the exemplary model problem of \Cref{ex:pde} with $g_{\mathrm{bc}}\equiv 0$. We emphasize, however, that similar constructions are possible for much more general problems.

As a first step, we transform the PDE~\eqref{eq:ellPDE} into a weak formulation by multiplying with an appropriate test function $w \in H^1_0(\Omega)$ and integrating by parts. We then seek a function $u \in H^1_0(\Omega)$ such that
\begin{equation}\label{eq:ellPDEweak}
    \la A \nabla u, \nabla w\ra_{L^2(\Omega)}
    =
    \la f, w \ra_{L^2(\Omega)},
    \qquad
    \text{for all } w\in H^1_0(\Omega).
\end{equation}
For convenience, we introduce the energy inner product
\begin{equation}\label{eq:bilinear-form-a}
    a(u,w)
    \coloneqq
    \la A \nabla u, \nabla w\ra_{L^2(\Omega)}.
\end{equation}
The weak formulation~\eqref{eq:ellPDEweak} can then be written as
\begin{equation}
    a(u,w)
    =
    \la f, w \ra_{L^2(\Omega)},
    \qquad
    \text{for all } w\in H^1_0(\Omega).
\end{equation}

An advantage of the weak formulation is that it reduces the regularity assumptions on the solution: second derivatives are not required and the first derivative needs to exist in a weak sense only. Nonetheless, \eqref{eq:ellPDEweak} cannot be solved directly, as it is still posed in the infinite-dimensional space $H^1_0(\Omega)$. Therefore, we solve~\eqref{eq:ellPDEweak} in a finite-dimensional subspace $V_m \subset H^1_0(\Omega)$ with $\dim V_m = m$, that is, we seek $u_m \in V_m$ such that
\begin{equation}\label{eq:Galerkin}
    a(u_m,w_m)
    =
    \la f, w_m \ra_{L^2(\Omega)},
    \qquad
    \text{for all } w_m\in V_m.
\end{equation}
This approach is called the \emph{Galerkin method}.

Since $V_m$ is finite-dimensional, the Galerkin approximation can itself be
viewed as a particular instance of a parametric model. Let
$(\varphi_i)_{1\leq i\leq m}$ be a basis of $V_m$ and consider the linear
parametrization
\begin{equation}
    \paramLinear\in\R^m
    \longmapsto
    v_\paramLinear
    \coloneqq
    \sum_{i=1}^m \paramLinear_i\varphi_i
    \in V_m .
\end{equation}
Its image is exactly $V_m$, and, since the parametrization is linear, its
tangent space is independent of the coefficient vector. That is, 
\begin{equation}
    T_\paramLinear\calM_\pmFun
    =
    V_m
    \qquad
    \text{for all } \paramLinear\in\R^m.
\end{equation}
Thus, Galerkin methods fit naturally into the parametric-model framework: they
correspond to a linear, finite-dimensional model whose tangent space coincides
everywhere with the approximation space itself. Writing the Galerkin solution
as
\begin{equation}
    u_m
    =
    \sum_{i=1}^{m} \paramLinear_i\varphi_i,
\end{equation}
the variational problem~\eqref{eq:Galerkin} therefore reduces to a
finite-dimensional linear system for the parameter vector
$\paramLinear=(\paramLinear_i)_{i=1}^m\in\R^m$.

Typical examples for $V_m$ are so-called \emph{finite element spaces} that consist of piece-wise polynomial and continuous functions, see, e.g., the works~\cite{Cia78,Bra07,BreS08} on finite element methods.

More generally, one may choose a test space $W_q \subset H^1_0(\Omega)$ distinct from the approximation space~$V_m$. The associated \emph{Petrov--Galerkin method} seeks $u_m \in V_m$ such that
\begin{equation}\label{eq:PetrovGalerkin}
    a(u_m,w_q)
    =
    \la f, w_q \ra_{L^2(\Omega)}
    \qquad
    \text{for all } w_q\in W_q,
\end{equation}
where the Galerkin method is recovered for the choice $W_q = V_m$. Note that this formulation allows for $\dim V_m=m \neq q =\dim W_q$. Classical square Petrov--Galerkin discretizations usually take $m=q$; when $m\neq q$, the discrete system is rectangular and can instead be interpreted in a least-squares or minimum-norm sense when appropriate. This will become important in \Cref{sec:anagram-green}.

\subsection{Neural networks}
\label{subsec:neural-networks}

Before introducing the concept of so-called physics-informed neural networks, we first provide a functional-analytic definition of a neural network that fits within the abstract framework of a \emph{parametric model}, which has been introduced in~\Cref{subsec:parametric-models} above.

\begin{definition}[multi-layer perceptron architecture]\label{def:architecture}
Based on the \emph{hyperparameters}
\begin{itemize}
  \item $L\in\N$, called the \emph{depth} of the network,
  \item $(l_i)_{0\le i\le L}\in\N^{L+1}$, called
  the \emph{widths} of the layers,
  \item $(\activation_i)_{1\le i\le L}\in\prod_{i=1}^{L}W_{\mathrm{loc}}^{1,\infty}(\R^{l_i};\R^{l_i})$, called the \emph{activation functions},
  \item $p\in\N$ with $1\le p\le \sum_{i=1}^L l_i (l_{i-1}+1)$, called the \emph{number of parameters},
  \item $\archChart\in C^1\left(\R^p;\prod_{i=1}^{L}\R^{l_i\times l_{i-1}}\times \prod_{i=1}^{L}\R^{l_i}\right)$,
  called the \emph{parameterization map},
\end{itemize}
a \emph{multi-layer perceptron (MLP) architecture} is a mapping $v \colon \R^p \to C^0(\Omega;\R^{l_L})$ for some subdomain $\Omega \subset \R^{l_0}$, defined by
\begin{equation}
    \paramPM \mapsto\circ_{k=1}^L \activation_k\big(\archChart_k(\paramPM)\bullet+\archChart_{L+k}(\paramPM)\big),
  \label{eq:architecture-map}
\end{equation}
where $\archChart_k(\paramPM) \in \R^{l_k \times l_{k-1}}$ and $\archChart_{L+k}(\paramPM) \in \R^{l_k}$ denote
the \emph{weight matrix} and \emph{bias vector} of the $k$th layer.
\end{definition}

In view of~\Cref{def:architecture}, we emphasize that activation functions are typically scalar-valued functions that are applied component-wise. That is, for a given $k\in\N$ and an activation function~$\activation\colon\R^k\to\R^k$, there exists $\sigma\colon\R\to\R$ such that for $x \in \R^k$, we have
\begin{equation*}
\activation(x) = \big(\sigma(x_i)\big)_{1\le i\le k}\in\R^k.
\end{equation*}
Therefore, it is common to abuse the notation and to refer to~$\activation$ through its scalar component function~$\sigma$.

Furthermore, we note that the parameterization map $\archChart$ encodes both
\emph{structural constraints} on the weights and biases, including parameter
sharing as in convolutional or tied-weight architectures, as well as the
distinction between \emph{learnable} and \emph{fixed} parameters.
In this paper, we will mainly view $\archChart$ as a convenient re-parameterization of a real vector of dimension
\begin{equation}
 p = \sum_{i=1}^L l_i (l_{i-1}+1)
\end{equation}
to the corresponding collection of~$(W_i)_{i=1}^L \in \prod_{i=1}^{L}\R^{l_i \times l_{i-1}}$ and~$(b_i)_{i=1}^{L}\in\prod_{i=1}^{L}\R^{l_i}$.
The MLP architecture of \Cref{def:architecture} therefore provides a particular instance of the parametric-model framework introduced in~\Cref{subsec:parametric-models}. Indeed, once the architecture (that is, the depth, layer widths, activation functions, and parameterization map $\archChart$) has been fixed, each parameter vector~$\paramPM\in\R^p$ determines a unique realized function $v_\paramPM$. The network can thus be viewed as a map
\begin{equation*}
    v:\R^p\to\calH,
    \qquad
    \paramPM\mapsto v_\paramPM,
\end{equation*}
provided that the realizations belong to the chosen Hilbert space $\calH$. In contrast with the linear parametrization associated with a Galerkin space, the dependence of $v_\paramPM$ on $\paramPM$ is generally nonlinear because of the successive compositions with the activation functions. Consequently, the set of realizations
\begin{equation*}
    \calM_\pmFun
    =
    \{v_\paramPM\mid\paramPM\in\R^p\}
    \subset\calH
\end{equation*}
is in general a nonlinear parametric manifold. Its tangent space at $v_\paramPM$ is generated by the parameter derivatives,
\begin{equation*}
    T_\paramPM\calM_\pmFun
    =
    \spn\{
        \partial_i v_\paramPM
        \mid
        1\leq i\leq p
    \},
\end{equation*}
and therefore generally depends on the current parameter value~$\paramPM$. This is the main distinction with the Galerkin parametrization above, whose tangent space is the fixed approximation space~$V_m$ at every parameter value.

\begin{remark}[Other architectures]
Note that we restrict our presentation to MLPs. Many other neural network architectures have been proposed, such as convolutional neural networks~\cite{lecunGradientbasedLearningApplied1998,goodfellowDeepLearning2016}, recurrent neural networks, variants thereof~\cite{elmanFindingStructureTime1990,hochreiterLongShorttermMemory1997,choPropertiesNeuralMachine2014}, and transformers~\cite{vaswaniAttentionAllYou2017}. From the perspective adopted in this work, these architectures can often be viewed as particular parameterizations obtained by imposing structural constraints on a general neural network model. Since our analysis only relies on the resulting parametric map~\eqref{eq:architecture-map}, focusing on MLPs entails no significant loss of generality.
\end{remark}

\subsection{Physics-informed neural networks}
\label{subsec:pinns}

Physics-informed neural networks (PINNs) provide a natural example of the
compound parametric models introduced in \Cref{subsec:parametric-models}. For
the boundary value problem~\eqref{eq:PDEstrong}, we introduce
\begin{equation*}
    \calY
    \coloneqq
    L^2(\Omega)\times L^2(\partial\Omega)
\end{equation*}
and define
\begin{equation}
    \pmComp^{D,B}
    \coloneqq
    (D,B)\circ v:
    \R^p\to\calY,
    \qquad
    \pmComp_\paramPM^{D,B}
    =
    \bigl(
        D[v_\paramPM],
        B[v_\paramPM]
    \bigr).
    \label{eq:pinn-compound-model}
\end{equation}
With the target
\begin{equation*}
    y^{D,B}
    \coloneqq
    \bigl(
        f,
        g_{\mathrm{bc}}
    \bigr),
\end{equation*}
the associated functional residual is
\begin{equation*}
    \resFun_\paramPM^{D,B}
    \coloneqq
    \pmComp_\paramPM^{D,B}-y^{D,B}.
\end{equation*}

The PINN formulation consists of minimizing its squared norm,
\begin{equation}\label{eq:lossPINNfunctional}
    \ell(\paramPM)
    =
    \norm[\resFun_\paramPM^{D,B}]_{\calY}^2
    =
    \int_\Omega (D[v_\paramPM](x)-f(x))^2 \,\dd x
    +
    \int_{\partial\Omega}(B[v_\paramPM](s)-g_{\mathrm{bc}}(s))^2 \,\dd s.
\end{equation}
In practice, this functional residual is discretized by pointwise evaluations.
Let $X_\Omega = \{x_i\}_{i=1}^{N_\Omega}\subset\Omega$ and
$X_{\partial\Omega} = \{s_i\}_{i=1}^{N_{\partial\Omega}}\subset\partial\Omega$
be families of interior and boundary collocation points. The resulting PINN loss is given by
\begin{equation}\label{eq:lossPINN}
    \ell_\mathrm{pinn}(\paramPM)
    =
    \frac1{N_\Omega}\sum_{i=1}^{N_\Omega}
    (D[v_\paramPM](x_i)-f(x_i))^2
    +
    \frac1{N_{\partial\Omega}}\sum_{i=1}^{N_{\partial\Omega}}
    (B[v_\paramPM](s_i)-g_{\mathrm{bc}}(s_i))^2.
\end{equation}
Thus, standard collocation replaces the functional residual by a finite
collection of pointwise measurements.
These pointwise evaluations must of course be well defined. In particular,
the neural-network realization must possess the regularity required by the
operators $D$ and $B$. This
typically requires sufficiently smooth activation functions. Automatic
differentiation~\cite{linnainmaaTaylorExpansionAccumulated1976,
baydinAutomaticDifferentiationMachine2018}
is then used to evaluate these quantities efficiently
\cite{raissiPhysicsinformedNeuralNetworks2019}. We emphasize that this
pointwise requirement is stronger than merely assuming
$\resFun_\paramPM^{D,B}\in\calY$; this distinction will become important when more
general measurements of the functional residual are introduced below.

For the specific problem~\eqref{eq:ellPDE} in \Cref{ex:pde} with
$g_{\mathrm{bc}}\equiv0$, the discrete loss becomes
\begin{equation}\label{eq:lossPINNell}
    \ell_\mathrm{pinn}(\paramPM)
    =
    \frac1{N_\Omega}\sum_{i=1}^{N_\Omega}
    (-\ddiv(A(x_i)\nabla v_\paramPM(x_i)) - f(x_i))^2
    +
    \frac1{N_{\partial\Omega}}\sum_{i=1}^{N_{\partial\Omega}}
    (\gamma v_\paramPM(s_i))^2,
\end{equation}
where $\gamma$ denotes the trace operator as before.

\section{Connecting Gauss--Newton with Petrov--Galerkin discretizations}
\label{sec:anagram-green}

We now introduce the second main ingredient of this work: the Gauss--Newton
(GN) method and its interpretation through Petrov--Galerkin discretizations.
We first recall the classical finite-dimensional construction. We then show
that a functional problem can be reduced to exactly the same setting by
choosing a finite family of measurements of the functional residual.
Different choices of approximation models and measurements recover, among
others, natural-gradient methods, kernel methods, regularization and
sketching strategies, both for PINNs and for direct regression problems.

\subsection{Gauss--Newton in finite-dimensional optimization}
\label{subsec:gauss-newton-finite-dimensional}

Let
\begin{equation*}
    \pmDis:
    \R^p
    \to
    \R^N,
    \qquad
    \paramPM
    \mapsto
    \pmDis_\paramPM,
\end{equation*}
be a differentiable discrete model and let $y\in\R^N$ be a target. We
consider the nonlinear least-squares problem
\begin{equation}
    \ell(\paramPM)
    =
    \frac12
    \norm[
        \pmDis_\paramPM-y
    ]_{\R^N}^2.
    \label{eq:gn-finite-loss}
\end{equation}
We denote its residual by
\begin{equation*}
    \resDis_\paramPM
    \coloneqq
    \pmDis_\paramPM-y
\end{equation*}
and the Jacobian of $\pmDis$ at $\paramPM$ by
\begin{equation*}
    \jacDis_\paramPM
    \coloneqq
    \dd\pmDis_\paramPM
    \in
    \R^{N\times p}.
\end{equation*}
Let $\bar\paramPM=\paramPM-\paramUpdate$, where $\paramUpdate\in\R^p$ denotes the correction to be
subtracted. The first-order approximation
\begin{equation}
    \resDis_{\paramPM-\paramUpdate}
    =
    \resDis_\paramPM
    -
    \jacDis_\paramPM\paramUpdate
    +
    o(\norm[\paramUpdate])
    \label{eq:gn-finite-linearization}
\end{equation}
leads to the local least-squares problem
\begin{equation}
    \paramUpdate_\paramPM^{\operatorname{GN}}
    \in
    \argmin_{\paramUpdate\in\R^p}
    \frac12
    \norm[
        \resDis_\paramPM-\jacDis_\paramPM\paramUpdate
    ]_{\R^N}^2.
    \label{eq:gn-finite-minimization}
\end{equation}

This formulation makes the basic idea of Gauss--Newton transparent. The
residual $\resDis_\paramPM$ is the correction that would ideally bring
$\pmDis_\paramPM$ to the target, whereas parameter variations can only
produce, as first-order approximations, corrections in
\begin{equation*}
    \Ima(\jacDis_\paramPM)
    =
    \spn\{
        \jacDis_\paramPM e_j
        \mid
        1\leq j\leq p
    \}.
\end{equation*}
Gauss--Newton therefore selects the accessible correction that best
approximates $\resDis_\paramPM$. Its optimality conditions are
\begin{equation}
    \left\langle
        \jacDis_\paramPM e_j,
        \resDis_\paramPM
        -
        \jacDis_\paramPM\paramUpdate_\paramPM^{\operatorname{GN}}
    \right\rangle_{\R^N}
    =
    0,
    \qquad
    1\leq j\leq p,
    \label{eq:gn-finite-testing}
\end{equation}
or, equivalently,
\begin{equation*}
    \jacDis_\paramPM^\top \jacDis_\paramPM
    \paramUpdate_\paramPM^{\operatorname{GN}}
    =
    \jacDis_\paramPM^\top \resDis_\paramPM.
\end{equation*}
Thus,
\begin{equation}
    \jacDis_\paramPM\paramUpdate_\paramPM^{\operatorname{GN}}
    =
    \Pi_{\Ima(\jacDis_\paramPM)}\resDis_\paramPM,
    \label{eq:gn-finite-projection}
\end{equation}
where $\Pi_{\Ima(\jacDis_\paramPM)}$ denotes the Euclidean orthogonal projection onto
the range of the Jacobian $\jacDis_\paramPM$. With the Moore--Penrose convention and denoting the pseudo-inverse with $(\bullet)^\dagger$, the
minimum-norm parameter correction is
\begin{equation}
    \paramUpdate_\paramPM^{\operatorname{GN}}
    =
    \jacDis_\paramPM^\dagger \resDis_\paramPM.
    \label{eq:gn-finite-pseudoinverse}
\end{equation}

\subsection{Discretized Gauss--Newton for functional models}
\label{subsec:discretized-functional-gn}

We now consider a parametric model as introduced in
\Cref{subsec:parametric-models},
\begin{equation*}
    \pmFun:
    \R^p
    \to
    \calH,
    \qquad
    \paramPM
    \mapsto
    \pmFun_\paramPM,
\end{equation*}
where $\calH$ is a Hilbert space, together with a target
$y\in\calH$. We denote the corresponding functional residual by
\begin{equation}
    \resFun_\paramPM
    \coloneqq
    \pmFun_\paramPM-y.
    \label{eq:functional-residual-general}
\end{equation}
For a parameter correction $\paramUpdate\in\R^p$,
\begin{equation*}
    \resFun_{\paramPM-\paramUpdate}
    =
    \resFun_\paramPM
    -
    \dd\pmFun_\paramPM[\paramUpdate]
    +
    o(\norm[\paramUpdate]).
\end{equation*}
At the linearized level, the ideal correction would therefore satisfy
\begin{equation}
    \dd\pmFun_\paramPM[\paramUpdate]
    =
    \resFun_\paramPM.
    \label{eq:functional-linearized-equation}
\end{equation}
The admissible corrections on the left-hand side belong to the
finite-dimensional tangent space
\begin{equation*}
    T_\paramPM\calM_\pmFun
    =
    \Ima(\dd\pmFun_\paramPM)
    =
    \spn\{
        \partial_j\pmFun_\paramPM
        \mid
        1\leq j\leq p
    \}
    \subset\calH.
\end{equation*}
Thus, although~\eqref{eq:functional-linearized-equation} is an equation
between functional objects, the space of admissible first-order corrections
is finite-dimensional.

To discretize this equation, we fix the current parameter $\paramPM$ and
introduce the \emph{extended tangent space}
\begin{equation}
    \overline{T_\paramPM\calM_\pmFun}
    \coloneqq
    \spn\{\resFun_\paramPM\}
    +
    T_\paramPM\calM_\pmFun
    \subset
    \calH.
    \label{eq:extended-tangent-space}
\end{equation}
Its dimension is at most $p+1$.
Let 
\begin{equation*}
    \Lambda_\paramPM
    =
    (\linform_i^\paramPM)_{i=1}^N,
    \qquad
    \linform_i^\paramPM
    \in
    \left(\overline{T_\paramPM\calM_\pmFun}\right)^{\prime},
\end{equation*}
be a family of linear measurements on
$\overline{T_\paramPM\calM_\pmFun}$.
Applying these measurements to
\eqref{eq:functional-linearized-equation} gives the finite-dimensional system
\begin{equation}
    \jacDis_\paramPM^{\Lambda_\paramPM}\paramUpdate
    =
    \resDis_\paramPM^{\Lambda_\paramPM},
    \qquad
    \resDis_{\paramPM,i}^{\Lambda_\paramPM}
    \coloneqq
    \linform_i^\paramPM(\resFun_\paramPM),
    \qquad
    \jacDis_{\paramPM,ij}^{\Lambda_\paramPM}
    \coloneqq
    \linform_i^\paramPM(
        \partial_j\pmFun_\paramPM
    ).
    \label{eq:discrete-functional-system}
\end{equation}
The linearized problem in the function space~$\calH$ has therefore been reduced to an ordinary
finite-dimensional problem of precisely the form considered in
\Cref{subsec:gauss-newton-finite-dimensional}. The corresponding
\emph{discretized Gauss--Newton problem} is defined by
\begin{equation}
    \paramUpdate_\paramPM^{\Lambda_\paramPM}
    \in
    \argmin_{\paramUpdate\in\R^p}
    \frac12
    \norm[
        \resDis_\paramPM^{\Lambda_\paramPM}
        -
        \jacDis_\paramPM^{\Lambda_\paramPM}\paramUpdate
    ]_{\R^N}^2.
    \label{eq:discretized-gn}
\end{equation}
As before, the minimum-norm correction is
\begin{equation}
    \paramUpdate_\paramPM^{\Lambda_\paramPM}
    =
    \bigl(
        \jacDis_\paramPM^{\Lambda_\paramPM}
    \bigr)^\dagger
    \resDis_\paramPM^{\Lambda_\paramPM}.
    \label{eq:discretized-gn-pseudoinverse}
\end{equation}
This construction has an immediate Petrov--Galerkin interpretation.
More generally, suppose that the measurements admit a representation
through a pairing with a test space $\calZ$, namely, for $1 \leq i  \leq N$,
\begin{equation}
    \linform_i^\paramPM(w)
    =
    \braket{
        w
    }{
        \testform_i^\paramPM
    }_{\calH\times\calZ},
    \qquad
    \text{for all } w\in\overline{T_\paramPM\calM_\pmFun},
    \label{eq:measurement-test-pairing}
\end{equation}
for some test function $\testform_i^\paramPM\in\calZ$.
Consequently, \eqref{eq:discrete-functional-system} can be written
equivalently as
\begin{equation}
    \braket{
        \dd\pmFun_\paramPM[\paramUpdate]
        -
        \resFun_\paramPM
    }{
        \testform_i^\paramPM
    }_{\calH\times\calZ}
    =
    0,
    \qquad
    1\leq i\leq N.
    \label{eq:functional-petrov-galerkin}
\end{equation}
Thus, at a fixed parameter $\paramPM\in\R^p$, choosing measurements of the
linearized problem in~$\calH$ is equivalent, whenever such a representation
is available, to choosing the Petrov--Galerkin test functions against which
this problem is enforced.

In particular, when $\calZ=\calH$ and the pairing is the inner product of
$\calH$, the finite dimensionality of
$\overline{T_\paramPM\calM_\pmFun}$ implies that every linear measurement is
continuous and admits a unique Riesz representer
$\rieszform_i^\paramPM\in\overline{T_\paramPM\calM_\pmFun}$, so that one may
take
\begin{equation*}
    \testform_i^\paramPM
    =
    \rieszform_i^\paramPM.
\end{equation*}
When the resulting system is rectangular or rank deficient, Gauss--Newton
provides the corresponding least-squares or minimum-norm solution through
\eqref{eq:discretized-gn}.

Finally, when the measurements extend to a common linear subspace
$\calD\subset\calH$ containing the image of $\pmFun$ and the target $y$, and which
are independent of $\paramPM$, they define the observation map
\begin{equation*}
    \calE_\Lambda:
    \calD
    \to
    \R^N,
    \qquad
    w
    \mapsto
    \bigl(
        \linform_i(w)
    \bigr)_{i=1}^N,
\end{equation*}
and hence a global discrete model
\begin{equation}
    \pmDis^\Lambda
    \coloneqq
    \calE_\Lambda\circ\pmFun:
    \R^p
    \to
    \R^N,
    \qquad
    y^\Lambda
    \coloneqq
    \calE_\Lambda(y).
    \label{eq:discrete-compound-model}
\end{equation}
The construction above is then exactly the ordinary Gauss--Newton
linearization of the discrete least-squares problem associated with
$\pmDis^\Lambda$ and the target $y^\Lambda$.
The local formulation is more general, though, since the measurements only need to be defined on the current extended tangent space.

\subsection{Recovering weak formulations}
\label{subsec:gn-weak-formulations}

We now show how classical weak formulations arise naturally within the
functional Gauss--Newton framework through suitable choices of linear
measurements and associated test functions. In particular, operator-adapted
Green sections and general weak test functions allow the functional residual
equation to be discretized without requiring a pointwise strong-form
residual.

More broadly, the same framework recovers a variety of classical
optimization and discretization methods by varying the parametric
approximation model and the measurement family. Examples including natural
gradient, kernel methods, empirical natural gradient, pointwise
Gauss--Newton, sketching, and regularization are briefly discussed in
\Cref{app:gn-approximation-and-test-examples}.

\subsubsection{Green sections.}

We consider the elliptic problem of \Cref{ex:pde} with homogeneous Dirichlet
boundary conditions, and denote with
\begin{equation}
    \wLapA:
    H_0^1(\Omega)
    \to
    H^{-1}(\Omega)
    \label{eq:weak-elliptic-operator}
\end{equation}
the weak elliptic operator defined by
\begin{equation*}
    \braket{
        \wLapA w
    }{
        z
    }_{H^{-1}(\Omega)\times H_0^1(\Omega)}
    =
    a(w,z).
\end{equation*}
Given a functional parametric model
\begin{equation*}
    \pmFun:
    \R^p
    \to
    H_0^1(\Omega),
\end{equation*}
we consider the compound model
\begin{equation}
    \pmComp^A
    \coloneqq
    \wLapA\circ\pmFun:
    \R^p
    \to
    H^{-1}(\Omega)
    \label{eq:elliptic-compound-model}
\end{equation}
with residual
\begin{equation*}
    \resFun_\paramPM^A
    =
    \wLapA\pmFun_\paramPM-f
    =
    \wLapA(\pmFun_\paramPM-u),
\end{equation*}
where $u$ denotes the exact solution to $\wLapA u=f$.
At the current parameter $\paramPM$, we consider the local solution space
\begin{equation*}
    \overline{T_\paramPM\calM_\pmFun}
    \coloneqq
    \spn\{\pmFun_\paramPM-u\}
    +
    T_\paramPM\calM_\pmFun
    \subset
    H_0^1(\Omega).
\end{equation*}
Assume that point evaluation at $x\in\Omega$ is well defined on this space.
Since $\overline{T_\paramPM\calM_\pmFun}$ is finite-dimensional, point evaluation
is then continuous and admits a corresponding local Green section
$\greenlift_{\paramPM,x}\in\overline{T_\paramPM\calM_\pmFun}$ satisfying
\begin{equation}
    a(w,\greenlift_{\paramPM,x})
    =
    w(x),
    \qquad
    \text{for all } w\in\overline{T_\paramPM\calM_\pmFun}.
    \label{eq:local-green-reproducing-property}
\end{equation}
Choosing $\greenlift_{\paramPM,x_i}$ as test functions in the Petrov--Galerkin
formulation gives
\begin{equation*}
    \braket{
        \dd\pmComp_\paramPM^A[\paramUpdate]
        -
        \resFun_\paramPM^A
    }{
        \greenlift_{\paramPM,x_i}
    }_{H^{-1}(\Omega)\times H_0^1(\Omega)}
    =
    0.
\end{equation*}
Since
\begin{equation*}
    \dd\pmComp_\paramPM^A[\paramUpdate]
    =
    \wLapA\dd\pmFun_\paramPM[\paramUpdate],
\end{equation*}
the reproducing property yields
\begin{equation}
    \dd\pmFun_\paramPM[\paramUpdate](x_i)
    =
    \pmFun_\paramPM(x_i)-u(x_i),
    \qquad
    1\leq i\leq N.
    \label{eq:green-interpolation-condition}
\end{equation}
Thus, Green sections transform weak measurements of the PDE residual into
pointwise measurements of the desired correction in the solution space.
Importantly, the differential operator is transferred from the approximation
model to the test functions, so that it does not need to be evaluated
directly on the model.

\subsubsection{General weak test functions.}
\label{subsubsec:general-weak-test-functions}

Green sections are only one particular choice of weak test functions.
We consider again the compound model
\begin{equation*}
    \pmComp^A
    \coloneqq
    \wLapA\circ\pmFun:
    \R^p
    \to
    H^{-1}(\Omega),
\end{equation*}
with functional residual
\begin{equation*}
    \resFun_\paramPM^A
    \coloneqq
    \pmComp_\paramPM^A-f
    =
    \wLapA\pmFun_\paramPM-f,
    \qquad
    f\in H^{-1}(\Omega).
\end{equation*}
Let $(\testform_i)_{i=1}^N$ be arbitrary functions in
$H_0^1(\Omega)$. At the current parameter $\paramPM$, they define the local
measurements
\begin{equation*}
    \linform_i^\paramPM(r)
    \coloneqq
    \braket{
        r
    }{
        \testform_i
    }_{H^{-1}(\Omega)\times H_0^1(\Omega)}
\end{equation*}
on the extended tangent space of $\pmComp^A$.
The discretized functional Gauss--Newton equation therefore becomes
\begin{equation}
    \braket{
        \dd\pmComp_\paramPM^A[\paramUpdate]
        -
        \resFun_\paramPM^A
    }{
        \testform_i
    }_{H^{-1}(\Omega)\times H_0^1(\Omega)}
    =
    0,
    \qquad
    1\leq i\leq N.
    \label{eq:weak-petrov-galerkin}
\end{equation}
Using the definition of $\wLapA$, this is equivalent to
\begin{equation}
    a(
        \dd\pmFun_\paramPM[\paramUpdate],
        \testform_i
    )
    =
    a(
        \pmFun_\paramPM,
        \testform_i
    )
    -
    \braket{
        f
    }{
        \testform_i
    }_{H^{-1}(\Omega)\times H_0^1(\Omega)},
    \qquad
    1\leq i\leq N.
    \label{eq:weak-petrov-galerkin-elliptic}
\end{equation}
When $f\in L^2(\Omega)$, the last duality pairing reduces to the usual
$L^2(\Omega)$ inner product.

The corresponding finite residual and Jacobian components in
\eqref{eq:discrete-functional-system} are therefore
\begin{equation}
    \resDis_{\paramPM,i}^{Z}
    \coloneqq
    a(
        \pmFun_\paramPM,
        \testform_i
    )
    -
    \braket{
        f
    }{
        \testform_i
    }_{H^{-1}(\Omega)\times H_0^1(\Omega)},
    \label{eq:weak-residual-components}
\end{equation}
and
\begin{equation}
    \jacDis_{\paramPM,ij}^{Z}
    \coloneqq
    a(
        \partial_j\pmFun_\paramPM,
        \testform_i
    ).
    \label{eq:weak-jacobian-components}
\end{equation}
The parameter correction is then obtained by applying ordinary
finite-dimensional Gauss--Newton to these weak residual components.

This formulation makes explicit the regularity advantage of weak
measurements. Although the compound residual takes values in
$H^{-1}(\Omega)$, it is tested against functions in $H_0^1(\Omega)$.
For a second-order elliptic operator, the bilinear form
$a(\cdot,\cdot)$ involves only weak spatial derivatives of first order, so that the
strong differential operator never needs to be explicitly evaluated in the
approximation model. In this sense, derivatives are transferred to the weak
form and the test functions, allowing residuals of substantially weaker
regularity than in a pointwise strong-form discretization.

The Green-section construction of the previous paragraph is recovered as a
particular operator-adapted choice of these weak test functions. By the
Lax--Milgram theorem, the weak elliptic operator
\begin{equation*}
    \wLapA:
    H_0^1(\Omega)
    \to
    H^{-1}(\Omega)
\end{equation*}
is an isomorphism and therefore admits a bounded inverse. Hence, for any
\begin{equation*}
    \lambda
    \in
    H^{-1}(\Omega),
\end{equation*}
we may define its Green lift
\begin{equation*}
    \greenlift_\lambda
    \coloneqq
    \wLapA^{-1}[\lambda]
    \in
    H_0^1(\Omega).
\end{equation*}
By construction, 
\begin{equation*}
    a(\greenlift_\lambda,w)
    =
    \braket{
        \lambda
    }{
        w
    }_{H^{-1}(\Omega)\times H_0^1(\Omega)}
    \qquad
    \text{for all } w\in H_0^1(\Omega).
\end{equation*}
Consequently, choosing
$\testform_i=\greenlift_{\lambda_i}$ in
\eqref{eq:weak-petrov-galerkin} transforms the weak residual measurements
into the corresponding measurements of the desired correction in the solution
space, namely
\begin{equation}
    \lambda_i\bigl(
        \dd\pmFun_\paramPM[\paramUpdate]
    \bigr)
    =
    \lambda_i\bigl(
        \pmFun_\paramPM-u
    \bigr),
    \qquad
    1\leq i\leq N.
    \label{eq:green-lift-measurement-correction}
\end{equation}
Arbitrary weak test functions therefore extend the same principle beyond point
evaluations while retaining the low-regularity $H^{-1}$ formulation.

The test functions can consequently be selected according to the variational
structure, regularity, or computational properties of the problem. In the
next section, we combine this freedom with a hybrid finite
element--neural approximation model to obtain a hybrid weak-test
method.

\section{Hybrid finite element--neural method}
\label{sec:hybrid-method}

The preceding framework shows that optimization and discretization methods
can be constructed by varying the approximation model and the test functions
in the common Petrov--Galerkin formulation
\eqref{eq:functional-petrov-galerkin}.
We now exploit this freedom to
construct a hybrid finite element--neural method in which a finite element component
accounts exactly for a prescribed finite-dimensional space, while the neural
component is restricted to its $a$-orthogonal complement.

\subsection{Projected hybrid approximation model}
\label{subsec:hybrid-mathematical-description}

Let
\begin{equation*}
    V_n
    =
    \spn\{\varphi_i\}_{i=1}^n
    \subset
    H_0^1(\Omega)
\end{equation*}
be a finite-dimensional approximation space. In the following, we primarily
consider finite element spaces of moderate dimension, although the
construction does not depend on this particular choice. We denote by
\begin{equation*}
    \frakP_n^a:
    H_0^1(\Omega)
    \to
    V_n
\end{equation*}
the $a$-orthogonal projection onto $V_n$, and set
\begin{equation*}
    \frakQ_n^a
    \coloneqq
    I-\frakP_n^a.
\end{equation*}
Hence,
\begin{equation*}
    \Ima(\frakQ_n^a)
    =
    V_n^{\perp,a}
    \coloneqq
    \left\{
        v\in H_0^1(\Omega)
        \mid
        a(v,w_n)=0
        \quad
        \text{for all } w_n\in V_n
    \right\}.
\end{equation*}

Let $u_n\in V_n$ denote the Galerkin approximation associated with $V_n$,
i.e.,
\begin{equation}
    a(u_n,w_n)
    =
    \braket{
        f
    }{
        w_n
    }_{H^{-1}(\Omega)\times H_0^1(\Omega)}
    \qquad
    \text{for all } w_n\in V_n.
    \label{eq:hybrid-fem-reference}
\end{equation}
Starting from the functional neural model
\begin{equation*}
    \pmFun:
    \R^p
    \to
    H_0^1(\Omega),
    \qquad
    \paramPM
    \mapsto
    v_\paramPM,
\end{equation*}
we define the projected hybrid model
\begin{equation}
    \pmFun_n^{\mathrm{hyb}}:
    \R^p
    \to
    H_0^1(\Omega),
    \qquad
    \paramPM
    \mapsto
    \pmFun_{\paramPM,n}
    \coloneqq
    u_n+\frakQ_n^a\pmFun_\paramPM.
    \label{eq:hybrid-projected-model}
\end{equation}
Thus, the finite element component accounts for the directions in $V_n$,
while the trainable contribution is restricted to
$V_n^{\perp,a}$.
Equivalently,
\begin{equation}
    \pmFun_{\paramPM,n}
    =
    \pmFun_\paramPM
    +
    z_{n,\paramPM},
    \qquad
    z_{n,\paramPM}
    \coloneqq
    u_n-\frakP_n^a\pmFun_\paramPM
    \in V_n.
    \label{eq:hybrid-model}
\end{equation}
The finite element term $z_{n,\paramPM}$ therefore acts as a
\emph{compensator}, rather than as an additional set of trainable parameters.
By construction,
\begin{equation}
    a(\pmFun_{\paramPM,n},w_n)
    =
    \braket{
        f
    }{
        w_n
    }_{H^{-1}(\Omega)\times H_0^1(\Omega)}
    \qquad
    \text{for all } w_n\in V_n.
    \label{eq:fem-compensator}
\end{equation}
For each $\paramPM$, computing the compensator amounts to solving a
finite-dimensional Galerkin system whose matrix depends only on $V_n$ and
$a$, and can therefore be factorized once and reused throughout training.
The differential of the hybrid model with respect to $\paramPM$ is
\begin{equation}
    \dd\bigl(\pmFun_n^{\mathrm{hyb}}\bigr)_\paramPM
    =
    \frakQ_n^a\dd\pmFun_\paramPM,
    \label{eq:hybrid-projected-tangent}
\end{equation}
so that its tangent directions are themselves restricted to
$V_n^{\perp,a}$.
For the elliptic problem, we introduce the associated compound model
\begin{equation}
    \pmComp_n^{A,\mathrm{hyb}}
    \coloneqq
    \wLapA\circ\pmFun_n^{\mathrm{hyb}}:
    \R^p
    \to
    H^{-1}(\Omega),
    \label{eq:hybrid-compound-model}
\end{equation}
with residual
\begin{equation}
    r_{\paramPM,n}^A
    \coloneqq
    \pmComp_{\paramPM,n}^{A,\mathrm{hyb}}-f
    =
    \wLapA \pmFun_{\paramPM,n}-f.
    \label{eq:hybrid-residual}
\end{equation}
Equation~\eqref{eq:fem-compensator} immediately implies
\begin{equation}
    \braket{
        r_{\paramPM,n}^A
    }{
        w_n
    }_{H^{-1}(\Omega)\times H_0^1(\Omega)}
    =
    0
    \qquad
    \text{for all } w_n\in V_n.
    \label{eq:hybrid-residual-vh-orthogonal}
\end{equation}
The weak equation is therefore satisfied exactly on $V_n$, independently of
the neural parameters.

\subsection{Projected Petrov--Galerkin update}
\label{subsec:hybrid-projected-pg}

Let $(\testform_i)_{i=1}^N$ be arbitrary test functions in
$H_0^1(\Omega)$. Inserting the compound model
$\pmComp_n^{A,\mathrm{hyb}}$ into
\eqref{eq:functional-petrov-galerkin} gives
\begin{equation}
    \braket{
        \dd\bigl(\pmComp_n^{A,\mathrm{hyb}}\bigr)_\paramPM[\paramUpdate]
        -
        r_{\paramPM,n}^A
    }{
        \testform_i
    }_{H^{-1}(\Omega)\times H_0^1(\Omega)}
    =
    0,
    \qquad
    1\leq i\leq N,
    \label{eq:hybrid-linearized-pg-general}
\end{equation}
where $\paramUpdate\in\R^p$ denotes, as before, the correction subtracted from the
current parameter.

Both the residual and the tangent directions vanish when tested against
$V_n$. Consequently, only the complementary components of the test
functions matter, and we may replace each $\testform_i$ by
\begin{equation}
    \widetilde\testform_i
    \coloneqq
    \frakQ_n^a\testform_i
    \in
    V_n^{\perp,a}.
    \label{eq:projected-test-function}
\end{equation}
Since $\widetilde\testform_i\in V_n^{\perp,a}$,
\begin{equation*}
    a(
        \frakQ_n^a\dd\pmFun_\paramPM[\paramUpdate],
        \widetilde\testform_i
    )
    =
    a(
        \dd\pmFun_\paramPM[\paramUpdate],
        \widetilde\testform_i
    ),
\end{equation*}
and, since $u_n\in V_n$,
\begin{equation*}
    \braket{
        r_{\paramPM,n}^A
    }{
        \widetilde\testform_i
    }_{H^{-1}(\Omega)\times H_0^1(\Omega)}
    =
    a(
        \pmFun_\paramPM,
        \widetilde\testform_i
    )
    -
    \braket{
        f
    }{
        \widetilde\testform_i
    }_{H^{-1}(\Omega)\times H_0^1(\Omega)}.
\end{equation*}
Therefore, \eqref{eq:hybrid-linearized-pg-general} reduces to
\begin{equation}
    a(
        \dd\pmFun_\paramPM[\paramUpdate],
        \widetilde\testform_i
    )
    =
    a(
        \pmFun_\paramPM,
        \widetilde\testform_i
    )
    -
    \braket{
        f
    }{
        \widetilde\testform_i
    }_{H^{-1}(\Omega)\times H_0^1(\Omega)},
    \qquad
    1\leq i\leq N.
    \label{eq:hybrid-complement-equation}
\end{equation}
Thus, the proposed method is again precisely of the form of
\eqref{eq:functional-petrov-galerkin}, now with the approximation model
$\pmComp_n^{A,\mathrm{hyb}}$ and projected test functions
$\widetilde\testform_i=\frakQ_n^a\testform_i$.

In the notation of~\eqref{eq:discrete-functional-system}, the corresponding
finite residual and Jacobian are
\begin{equation}
    \resDis_{\paramPM,i}^n
    \coloneqq
    a(
        \pmFun_\paramPM,
        \widetilde\testform_i
    )
    -
    \braket{
        f
    }{
        \widetilde\testform_i
    }_{H^{-1}(\Omega)\times H_0^1(\Omega)},
    \label{eq:hybrid-feature-residual}
\end{equation}
and
\begin{equation}
    \jacDis_{\paramPM,ij}^n
    \coloneqq
    a(
        \partial_j\pmFun_\paramPM,
        \widetilde\testform_i
    ).
    \label{eq:hybrid-feature-jacobian}
\end{equation}
The parameter correction is therefore obtained from a similar
finite-dimensional Gauss--Newton problem as before and reads
\begin{equation}
    \jacDis_\paramPM^n\paramUpdate
    =
    \resDis_\paramPM^n,
    \label{eq:hybrid-linear-system}
\end{equation}
understood in the least-squares sense when the system is rectangular or
rank deficient.

The whole construction can be summarized as follows. The finite element
component is treated exactly in $V_n$, while
$\frakQ_n^a=I-\frakP_n^a$ removes the corresponding finite element directions
from both the neural correction and the test functions. Consequently, the
remaining optimization problem is posed entirely on the $a$-orthogonal
complement $V_n^{\perp,a}$.

\subsection{Numerical realization}
\label{subsec:hybrid-numerical-realization}

The hybrid formulation does not require an alternating optimization
between finite element and neural variables. The Galerkin solution $u_n$ and
a factorization of the finite element stiffness matrix are computed once.
For a given parameter $\paramPM$, evaluating the hybrid model only requires
computing the projection $\frakP_n^a\pmFun_\paramPM$, while projecting a test
function through~$\frakQ_n^a$ uses the same prefactorized finite-dimensional
system.

At each neural iteration, a family of test functions is selected and
projected onto $V_n^{\perp,a}$. The residual vector $\resDis_\paramPM^n$ and Jacobian
$\jacDis_\paramPM^n$ are then assembled from
\eqref{eq:hybrid-feature-residual} and \eqref{eq:hybrid-feature-jacobian}, respectively, defining an ordinary finite-dimensional Gauss--Newton problem.

After updating the neural parameters, the new approximation is simply
evaluated through
\begin{equation*}
    \pmFun_{\paramPM,n}
    =
    u_n+\frakQ_n^a\pmFun_\paramPM.
\end{equation*}
Recomputing the projection at the updated parameter is therefore part of
evaluating the hybrid approximation model itself, and should not be
interpreted as a block-alternating optimization step.

\section{Numerical experiments}
\label{sec:numerical-experiments}

The experiments are organized around two questions. First, we verify that the
functional discretization developed above allows Gauss--Newton algorithms
designed for finite residual vectors to be applied directly to weak PDE
formulations. Second, we examine whether this construction remains effective
beyond Green sections and smooth solutions, and whether the hybrid finite element--neural model improves the robustness of weak neural solvers.

Throughout the section, we monitor both relative $L^2$ and full relative
$H^1$ errors. Complementary implementation details, optimizer and activation
ablations, additional figures, and complete method-level results are provided
in
\weblink{https://nilo.schwencke.me/tutorials/beyond-pinns-companion/}{the companion blog}. 

\subsection{Gauss--Newton optimization of weak residuals}
\label{subsec:experiments-energy-vs-residual}

\subsubsection{Benchmark problems and weak-residual discretization}
We first consider two smooth one-dimensional problems on
$\Omega=(-1,1)$ with homogeneous Neumann boundary conditions.
Both may be written as
\begin{equation}
    -u''+\rho(u)=f,
    \qquad
    u'(-1)=u'(1)=0,
    \label{eq:exp-smooth-problems}
\end{equation}
with
\begin{equation*}
    \rho(s)=s
    \qquad\text{or}\qquad
    \rho(s)=s^3,
\end{equation*}
corresponding respectively to the linear and nonlinear benchmarks. In both
cases, the solution is
\begin{equation*}
    u(x)=\cos(\pi x),
\end{equation*}
so that
\begin{equation*}
    f(x)
    =
    \pi^2\cos(\pi x)
    +
    \rho\bigl(\cos(\pi x)\bigr).
\end{equation*}

The natural energy space is $H^1(-1,1)$, equipped with
\begin{equation*}
    a(v,w)
    =
    \int_{-1}^{1}
    \bigl(
        v'(x)w'(x)+v(x)w(x)
    \bigr)\,\dd x.
\end{equation*}
For $t\in[-1,1]$, let $g_t$ denote the Riesz representer of point
evaluation at $t$ with respect to $a$, given by
\begin{equation}
    g_t(x)
    =
    \frac{
        \cosh(\min\{x,t\}+1)
        \cosh(1-\max\{x,t\})
    }{\sinh(2)}.
    \label{eq:exp-smooth-green}
\end{equation}
For $t\in(-1,1)$, this is the Neumann Green section associated with
$-\partial_{xx}+1$.
Testing the residual against these functions and
integrating by parts, while retaining the Neumann conditions as separate
residual coordinates, gives
\begin{equation}
\begin{split}
    r_t^\rho(v)
    ={}&
    v(t)
    +
    v'(-1)g_t(-1)
    -
    v'(1)g_t(1)
    \\
    &+
    \int_{-1}^{1}
        \bigl(\rho(v(x))-v(x)\bigr)g_t(x)\,\dd x
    -
    \int_{-1}^{1}
        f(x)g_t(x)\,\dd x .
\end{split}
\label{eq:exp-smooth-weak-residual}
\end{equation}
For the linear problem, the first integral in the second line vanishes,
whereas for the nonlinear problem it contains $v^3-v$.

We discretize the weak PDE residual with $200$ fixed sections $g_t$ centered
at equispaced points of $[-1,1]$, including both endpoints, and append the two
Neumann residuals $v'(-1)$ and $v'(1)$, resulting in a residual vector with $202$ scalar components.
For interior centers $t$, $g_t'$ has a jump at $x=t$.
The integrals in~\eqref{eq:exp-smooth-weak-residual} are therefore split at $t$ and evaluated with a
$64$-point Gauss--Legendre rule on each subinterval.
These quadrature points are used only to evaluate the corresponding linear
measurements entering the residual vector and its Jacobian; they do not define
additional measurements and therefore do not increase the dimension of either.

\subsubsection{Energy references}
We compare ordinary Gauss--Newton (GN) solvers applied to this weak residual with
specialized energy-based methods. For the linear problem, the reference is the
Gauss--Newton Deep Ritz method of
\cite{haoGaussNewtonVariational2023}, associated with
\begin{equation}
    \mathcal E_{\mathrm{lin}}(v)
    =
    \frac12
    \int_{-1}^{1}
        \bigl(|v'(x)|^2+v(x)^2\bigr)\,\dd x
    -
    \int_{-1}^{1}f(x)v(x)\,\dd x .
    \label{eq:exp-linear-energy}
\end{equation}
Following the reference implementation, its integrals are evaluated by
two-point Gauss--Legendre quadrature on uniform cells with
$h_{\mathrm{train}}=1/3000$ and $h_{\mathrm{test}}=1/4000$, corresponding to
$12000$ and $16000$ quadrature nodes, respectively, and we use the same
geometric step-length search.

For the nonlinear problem, the reference is the Energy Natural Gradient
(NG) method of \cite{mullerAchievingHighAccuracy2023}, based on
\begin{equation}
    \mathcal E_{\mathrm{nl}}(v)
    =
    \frac12
    \int_{-1}^{1}|v'(x)|^2\,\dd x
    +
    \frac14
    \int_{-1}^{1}v(x)^4\,\dd x
    -
    \int_{-1}^{1}f(x)v(x)\,\dd x,
    \label{eq:exp-nonlinear-energy}
\end{equation}
whose second variation,
\begin{equation}
    D^2\mathcal E_{\mathrm{nl}}(v)[w,z]
    =
    \int_{-1}^{1}w'(x)z'(x)\,\dd x
    +
    3
    \int_{-1}^{1}v(x)^2w(x)z(x)\,\dd x,
    \label{eq:exp-nonlinear-energy-metric}
\end{equation}
defines the corresponding state-dependent natural-gradient metric. We report
both the published-protocol reference, using a width-$32$ $\tanh$ network and
trapezoidal integration with $20000$ training and $200000$ evaluation points,
and a matched version using the same metric with the controlled width-$64$
protocol.

\subsubsection{Residual Gauss--Newton solvers}
Let $\mathbf r_\theta\in\R^N$ denote a strong or weak residual vector and
let
\begin{equation*}
    J_\theta
    =
    U\Sigma V^\top,
    \qquad
    \Sigma
    =
    \operatorname{diag}(\sigma_1,\ldots,\sigma_q),
    \qquad
    q=\min\{N,p\},
\end{equation*}
be a thin singular-value decomposition of its Jacobian, with singular values in
non-increasing order. With the convention that the parameter update is
$\bar\theta=\theta-\xi$, we consider the cutoff and ridge Gauss--Newton
corrections
\begin{align}
    \xi_{\mathrm{cut}}
    &=
    V\,
    \operatorname{diag}
    \left(
        \frac{\mathbf 1_{\{\sigma_i\geq\tau\}}}{\sigma_i}
    \right)
    U^\top\mathbf r_\theta,
    \label{eq:exp-cutoff-gn}
    \\
    \xi_{\mathrm{ridge}}
    &=
    V\,
    \operatorname{diag}
    \left(
        \frac{\sigma_i}{\sigma_i^2+\lambda}
    \right)
    U^\top\mathbf r_\theta.
    \label{eq:exp-ridge-gn}
\end{align}
For ridge Gauss--Newton, the regularization parameter is held fixed throughout
each run. We use $\lambda=10^{-12}$ for the strong residuals, and
$\lambda=10^{-13}$ and $10^{-14}$ for the weak residuals in the linear and
nonlinear problems, respectively.
We also apply AMStramGRAM~\cite{schwencke2025amstramgram} and
DSGNAR~\cite{webbOptimisationFrameworkWellConditioned2026} directly to the
same finite residual vectors. No modification of these algorithms is required
for the weak formulation: once the test functionals have been evaluated, the
optimizer receives an ordinary residual vector together with its parameter
Jacobian.

For the smooth benchmarks, we additionally consider the corresponding
strong-residual formulations as controls.
Since the MLP trial functions use smooth $\tanh$ activations,
their pointwise second-order derivatives are well defined, allowing us to distinguish
the effect of the nonlinear least-squares geometry from that of weak testing.

\subsubsection{Experimental protocol and results}
All matched experiments use the same width-$64$ $\tanh$ architecture and are
repeated over ten seeds, with identical initial neural parameters across
methods for each seed.
The strong formulation uses $200$ pointwise PDE residuals and the weak
formulation $200$ Green-section residuals, evaluated at the same equispaced
points of $[-1,1]$;
the two
Neumann boundary residuals are appended in both cases.
Relative $L^2$ and full relative $H^1$ errors with respect to the solution
$u$ are evaluated using an independent high-resolution quadrature.
We report medians over the ten seeds,
with interquartile bands in the convergence plots; reported wall-clock times
exclude compilation and warm-up.

The principal weak-residual comparisons use ridge Gauss--Newton and DSGNAR.
Final errors and wall-clock times for the corresponding strong-residual
controls are also reported below. Additional results for cutoff
Gauss--Newton and AMStramGRAM, together with complete convergence curves and
further optimizer diagnostics, are provided in
\weblink{https://nilo.schwencke.me/tutorials/beyond-pinns-companion/}{the companion blog}.

\begin{figure}[t]
    \centering
    \includegraphics[width=\textwidth]
    {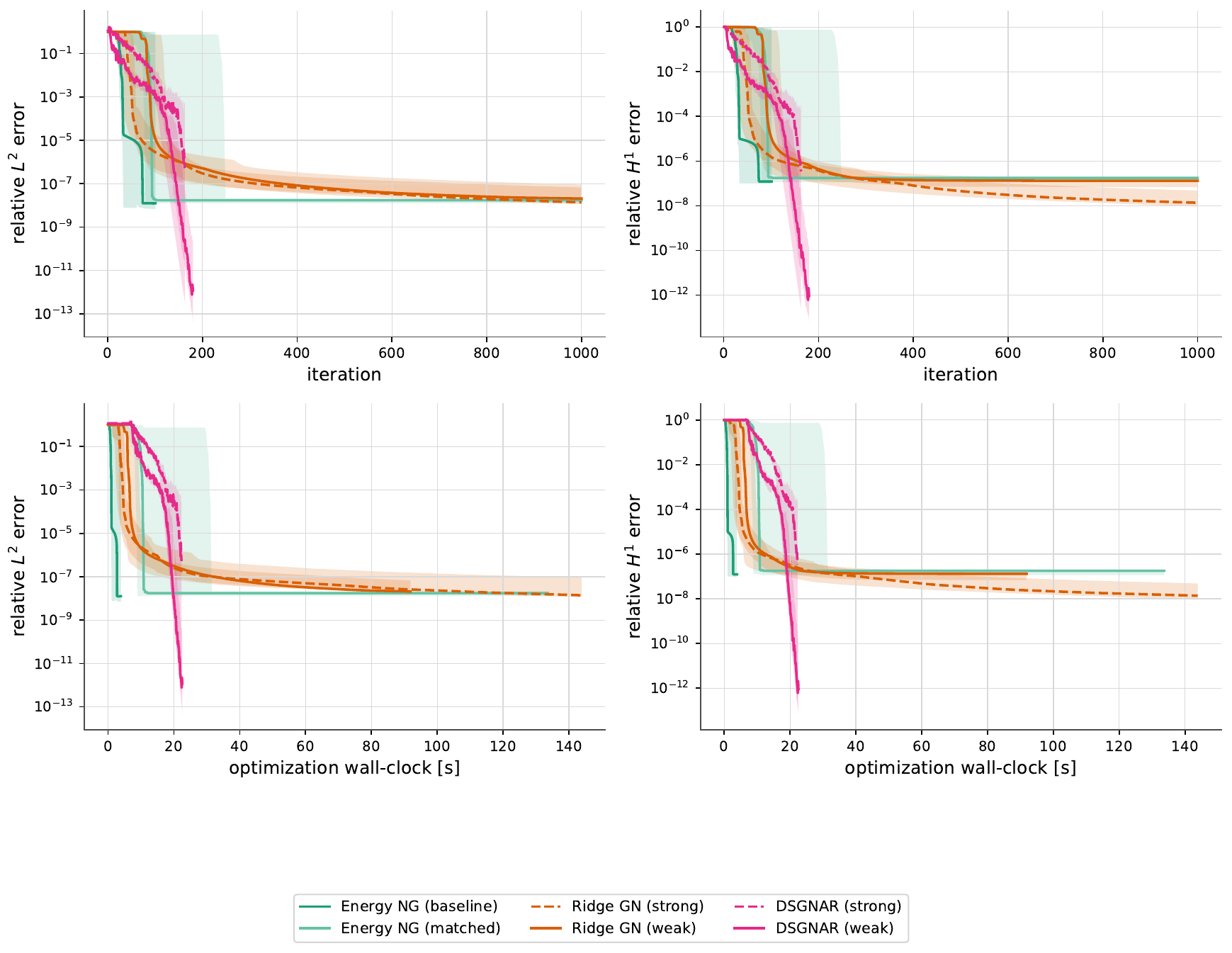}
    \caption{%
        \textbf{Nonlinear smooth benchmark with $\tanh$ networks.}
        Median relative $L^2$ and full relative $H^1$ errors over ten seeds,
        with interquartile bands, against iterations and optimization
        wall-clock. Energy NG (baseline) follows the protocol of
        \cite{mullerAchievingHighAccuracy2023}, whereas Energy NG (matched)
        uses the controlled width-$64$ protocol.
    }
    \label{fig:exp-smooth-nonlinear-formulations}
\end{figure}

The linear benchmark exhibits essentially the same qualitative behavior as
the nonlinear one; its complete convergence curves are provided in
\weblink{https://nilo.schwencke.me/tutorials/beyond-pinns-companion/}{the companion blog}.
We therefore show only the nonlinear convergence curves in
\Cref{fig:exp-smooth-nonlinear-formulations}, while
\Cref{tab:exp-smooth-summary} reports the final results for both problems.

\begin{table}[t]
    \centering
    \small
    \caption{%
        Median final relative errors and optimization wall-clock times over
        ten seeds for the smooth $\tanh$ benchmarks. Energy NG (baseline)
        follows the published protocol of
        \cite{mullerAchievingHighAccuracy2023}, whereas Energy NG (matched)
        uses the controlled width-$64$ protocol.
        Wall-clock times reflect the stated protocol of each method and are
        therefore not all intended as matched computational-cost comparisons.
    }
    \label{tab:exp-smooth-summary}
    \begin{tabular}{llccc}
        \toprule
        Problem
        & Method
        & rel.~$L^2$
        & rel.~$H^1$
        & time [s] \\
        \midrule
        \multirow{5}{*}{Linear}
        & GN Deep Ritz (baseline)
        & $1.12\times10^{-8}$
        & $1.04\times10^{-7}$
        & $127.3$ \\
        & ridge GN strong
        & $5.61\times10^{-8}$
        & $3.89\times10^{-8}$
        & $140.5$ \\
        & ridge GN weak
        & $6.16\times10^{-8}$
        & $5.15\times10^{-7}$
        & $76.9$ \\
        & DSGNAR strong
        & $3.22\times10^{-15}$
        & $1.36\times10^{-15}$
        & $24.6$ \\
        & DSGNAR weak
        & $6.70\times10^{-15}$
        & $4.33\times10^{-14}$
        & $21.0$ \\
        \midrule
        \multirow{6}{*}{Nonlinear}
        & Energy NG (baseline)
        & $1.27\times10^{-8}$
        & $1.23\times10^{-7}$
        & $3.9$ \\
        & Energy NG (matched)
        & $1.74\times10^{-8}$
        & $1.77\times10^{-7}$
        & $133.6$ \\
        & ridge GN strong
        & $1.39\times10^{-8}$
        & $1.37\times10^{-8}$
        & $143.8$ \\
        & ridge GN weak
        & $2.07\times10^{-8}$
        & $1.31\times10^{-7}$
        & $91.8$ \\
        & DSGNAR strong
        & $6.30\times10^{-15}$
        & $5.46\times10^{-15}$
        & $28.5$ \\
        & DSGNAR weak
        & $3.31\times10^{-15}$
        & $2.81\times10^{-14}$
        & $24.3$ \\
        \bottomrule
    \end{tabular}
\end{table}

\Cref{tab:exp-smooth-summary} shows first that no specialized energy metric is
required to optimize the weak formulation. Ridge Gauss--Newton applied
directly to the weak residual reaches an accuracy comparable to the dedicated
energy methods
\cite{haoGaussNewtonVariational2023,mullerAchievingHighAccuracy2023},
whereas DSGNAR drives both weak problems close to machine precision.

The strong-residual controls reach comparable or better accuracies, with
DSGNAR reaching near-machine precision in both formulations.
Hence, on these smooth problems, the results indicate that the dominant
improvement is associated with the residual least-squares formulation and its
Gauss--Newton treatment, rather than with weak testing itself.
The significance of the weak formulation is that this same
optimization structure survives after integration by parts, without requiring
the strong differential operator to be evaluated pointwise on the
approximation model.

The nonlinear Energy-NG baseline additionally exhibits substantial
seed dependence, as visible in
\Cref{fig:exp-smooth-nonlinear-formulations}. Six of the ten runs reach the
high-accuracy regime, with final relative $H^1$ errors between approximately
$8\times10^{-8}$ and $1.3\times10^{-7}$, whereas four remain close to unit
relative error. Its median is therefore accurate while its interquartile band
remains broad. Under the matched protocol, all ten runs reach the accurate
branch.

\subsection{General weak formulations and hybrid finite element--neural models}
\label{subsec:experiments-petrov-hybrid}

We next investigate whether the framework remains effective when
operator-adapted Green sections are unavailable or suboptimal, and when the
right-hand side or solution has limited regularity. We consider four
benchmarks.

\emph{Multiscale diffusion (MS).}
On $\Omega=(0,1)$, we consider the one-dimensional instance of the elliptic
problem of \Cref{ex:pde},
\begin{equation}
    -\frac{\dd}{\dd x}
    \left(
        A_\varepsilon(x)\frac{\dd}{\dd x}u(x)
    \right)
    =
    1,
    \qquad
    u(0)=u(1)=0,
    \label{eq:exp-ms-problem}
\end{equation}
with an oscillating diffusion coefficient
\begin{equation}
    A_\varepsilon(x)
    =
    \frac{1}{2+\cos(2\pi x/\varepsilon)},
    \qquad
    \varepsilon=\frac1{16}.
    \label{eq:exp-ms-coefficient}
\end{equation}
The corresponding solution involves $\varepsilon$-dependent oscillations; see
also~\cite[Ch.~2]{MalP21}.

\emph{Jump forcing (JF).}
On $\Omega=(0,1)^2$, we consider the homogeneous Dirichlet Poisson problem
with manufactured solution
\begin{equation*}
    u(x,y)=p(x)y(1-y),
    \qquad
    p(x)
    =
    \begin{cases}
        \displaystyle
        \frac{x^2}{2}-\frac{3x}{8},
        & x<\frac12,
        \\[1ex]
        \displaystyle
        \frac{x-1}{8},
        & x\geq\frac12.
    \end{cases}
\end{equation*}
The corresponding right-hand side $f=-\Delta u$ is discontinuous across
$x=1/2$, while the solution has no distributional singularity.

\emph{Line source (LS).}
Again on $\Omega=(0,1)^2$, we prescribe
\begin{equation}
    u(x,y)
    =
    \min(x,1-x)y(1-y)
    +
    \frac12\sin(2\pi x)\sin(2\pi y),
    \label{eq:exp-ls-solution}
\end{equation}
which solves
\begin{equation}
    -\Delta u
    =
    f_{\mathrm{reg}}
    +
    2y(1-y)\,\delta_{\{x=1/2\}},
    \qquad
    f_{\mathrm{reg}}
    =
    2\min(x,1-x)
    +
    4\pi^2\sin(2\pi x)\sin(2\pi y).
    \label{eq:exp-ls-problem}
\end{equation}
Here
\begin{equation*}
    \left\langle
        \delta_{\{x=1/2\}},z
    \right\rangle
    =
    \int_0^1 z\!\left(\frac12,y\right)\,\dd y,
\end{equation*}
so the singular source is treated directly as a one-dimensional functional
and is never replaced by a regularized volumetric source.

\emph{Reentrant corner (RC).}
Finally, on the L-shaped domain
\begin{equation}
    \Omega_L
    =
    (-1,1)^2
    \setminus
    \bigl([0,1)\times(-1,0]\bigr),
    \label{eq:exp-rc-domain}
\end{equation}
we consider the homogeneous Dirichlet Poisson problem with manufactured
solution
\begin{equation}
    u(x,y)
    =
    (1-x^2)(1-y^2)
    r^{2/3}
    \sin\!\left(\frac{2\vartheta}{3}\right),
    \label{eq:exp-rc-solution}
\end{equation}
where $(r,\vartheta)$ are polar coordinates centered at the reentrant corner.
The factor $r^{2/3}\sin(2\vartheta/3)$ gives the classical corner
singularity, with $|\nabla u|\sim r^{-1/3}$ near the origin.

\subsubsection{Trial and test spaces}

For MS, the pure neural model contains $1255$ trainable parameters and the
hybrid model combines a $375$-parameter neural component with $15$
interior $P_1$ finite element degrees of freedom, giving a total
approximation dimension of $390$. For JF and LS, the corresponding dimensions
are $6561$ for the pure neural model and $1605+165=1770$ for the hybrid model.
For RC, the hybrid model combines the same $1605$-parameter neural component
with $403$ graded finite element degrees of freedom, giving a total
approximation dimension of $2008$. The neural parametrizations use smooth $\tanh$
activation functions; MS, JF, and LS additionally use fixed Fourier features.
The finite element compensation spaces are deliberately coarse.

\begin{figure}[t]
    \centering

    \begin{subfigure}[t]{0.48\textwidth}
        \centering
        \includegraphics[width=\textwidth]
        {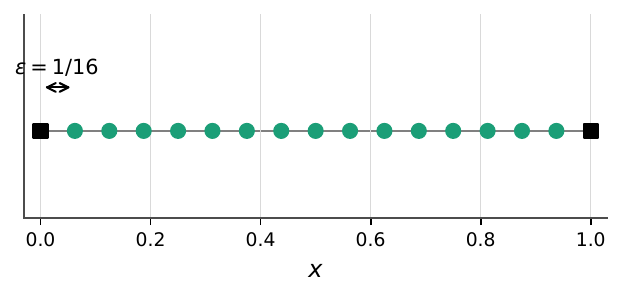}
        \caption{MS: $P_1$ mesh with $15$ interior degrees of freedom on
        $16$ elements.}
        \label{fig:exp-ms-compensation-mesh}
    \end{subfigure}
    \hfill
    \begin{subfigure}[t]{0.48\textwidth}
        \centering
        \includegraphics[width=\textwidth]
        {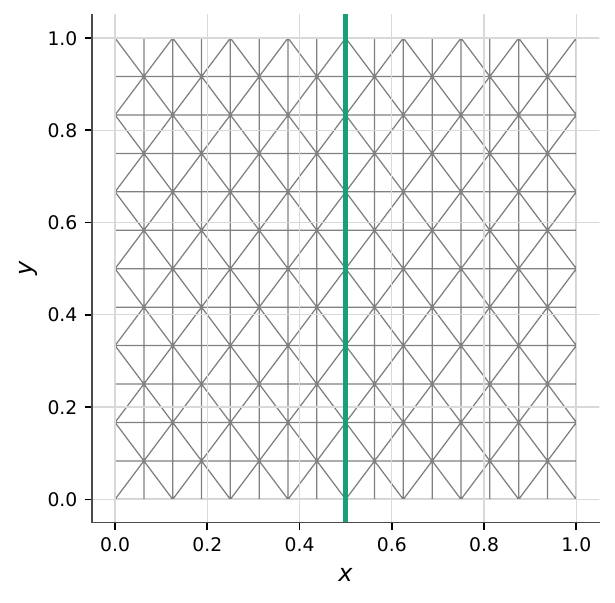}
        \caption{JF: interface-fitted triangular mesh, with $x=1/2$
        represented exactly.}
        \label{fig:exp-jf-compensation-mesh}
    \end{subfigure}

    \medskip

    \begin{subfigure}[t]{0.48\textwidth}
        \centering
        \includegraphics[width=\textwidth]
        {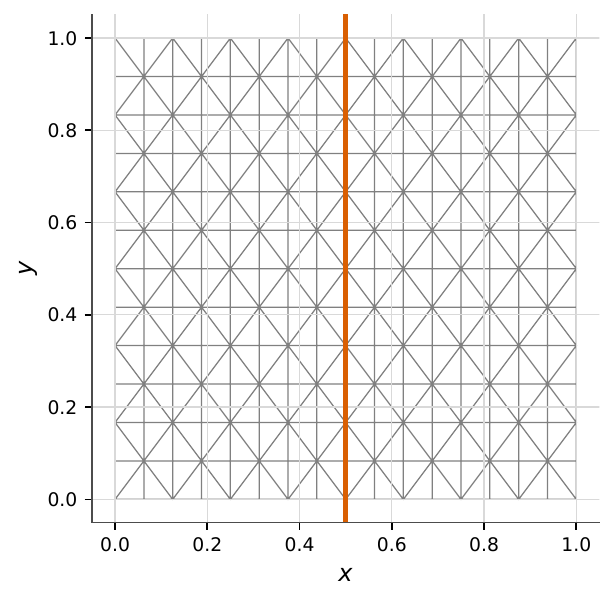}
        \caption{LS: interface-fitted triangular mesh, with the line source
        supported on $x=1/2$.}
        \label{fig:exp-ls-compensation-mesh}
    \end{subfigure}
    \hfill
    \begin{subfigure}[t]{0.48\textwidth}
        \centering
        \includegraphics[width=\textwidth]
        {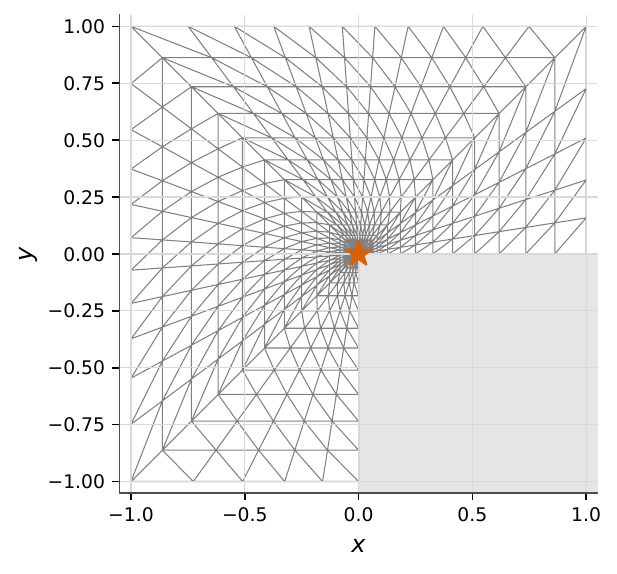}
        \caption{RC: graded $P_1$ mesh toward the reentrant corner.}
        \label{fig:exp-rc-compensation-mesh}
    \end{subfigure}

    \caption{%
        \textbf{Finite element compensation meshes.}
        The deliberately coarse finite element meshes used inside the hybrid
        models for MS, JF, LS, and RC.
        The JF and LS meshes resolve the line $x=1/2$ exactly, while the RC mesh is
        graded toward the reentrant corner.
        These meshes are distinct from the standalone finite element reference
        discretizations used in \Cref{tab:exp-general-summary}.
    }
    \label{fig:exp-petrov-compensation-meshes}
\end{figure}

The MS compensation space consists of $P_1$ functions on $16$
one-dimensional elements. For JF and LS, the compensation space is built on
a $16\times12$ interface-fitted triangular grid, so that $x=1/2$ is represented
exactly. For RC, the mesh is graded toward the reentrant corner.
The corresponding compensation meshes are shown in
\Cref{fig:exp-petrov-compensation-meshes}.

We compare several weak test families. In one dimension, besides the
operator-adapted $A_\varepsilon$-Green sections, we use the
$H_0^1(0,1)$ point-evaluation representers
\begin{equation}
    g_t(x)=\min(x,t)-xt.
    \label{eq:exp-h01-green}
\end{equation}
For JF and LS, if $(\lambda_k,\varphi_k)$ denote discrete
Dirichlet-Laplacian eigenpairs, we use truncated eigen-Green sections
\begin{equation}
    g_z^{(K)}(x)
    =
    \sum_{k=1}^K
    \frac{\varphi_k(z)\varphi_k(x)}{\lambda_k},
    \label{eq:exp-eigen-green}
\end{equation}
as well as tensor products of the one-dimensional functions
in~\eqref{eq:exp-h01-green}. For RC, we additionally consider the shifted
Sobolev family obtained by replacing $\lambda_k^{-1}$ in
\eqref{eq:exp-eigen-green} with $(\lambda_k+1)^{-1}$. In all four problems,
we also consider compact piecewise-affine hat functions with randomly chosen
locations.

The global test families use $4096$ test locations in one dimension and
$8192$ in two dimensions, while the local-hat experiments use $1000$ sampled
test functions. Fixed test functions are normalized in the corresponding energy norm.
The Fourier features belong only to the trial parametrization and are not
used as test functions.

Weak residuals are evaluated directly as variational functionals. For JF,
quadrature is split at the interface $x=1/2$ and, when necessary, at the
breakpoints of the test functions. For LS, the singular contribution is
evaluated separately as
\begin{equation*}
    2\int_0^1
        y(1-y)
        z\!\left(\frac12,y\right)
    \,\dd y.
\end{equation*}
For RC, both residual and error quadratures are graded toward the reentrant
corner. As in the smooth benchmarks, the quadrature points used to evaluate
a weak functional do not constitute additional residual measurements.

\subsubsection{Optimization and evaluation}

All trainable neural components are optimized with
DSGNAR~\cite{webbOptimisationFrameworkWellConditioned2026} in double
precision, with at most $300$ iterations over five seeds; we report median
relative $L^2$ and full relative $H^1$ errors evaluated by quadratures
independent of the training residual, with interquartile bands in convergence
plots.

Ablations comparing the exact projected hybrid formulation of
\Cref{sec:hybrid-method} with lagged-Jacobian and genuinely alternating
finite element--neural variants are provided in
\weblink{https://nilo.schwencke.me/tutorials/beyond-pinns-companion/}{the companion blog}.

\subsubsection{Results}
\label{subsubsec:hybrid-results}

\begin{figure}[t]
    \centering

    \begin{subfigure}[t]{0.48\textwidth}
        \centering
        \includegraphics[width=\textwidth]
        {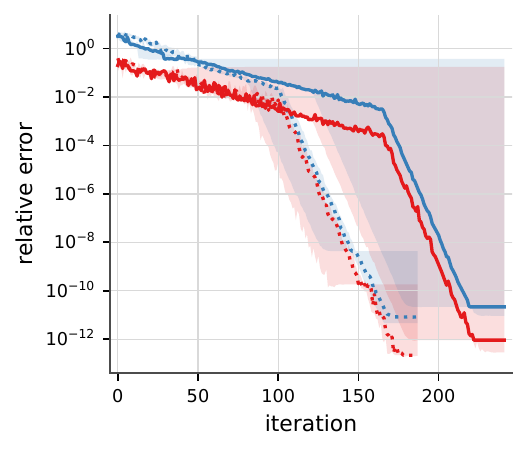}
        \caption{MS: multiscale diffusion.}
        \label{fig:exp-petrov-projection-ms}
    \end{subfigure}
    \hfill
    \begin{subfigure}[t]{0.48\textwidth}
        \centering
        \includegraphics[width=\textwidth]
        {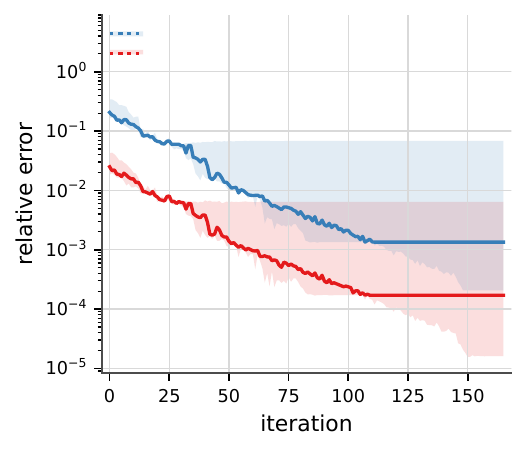}
        \caption{JF: jump forcing.}
        \label{fig:exp-petrov-projection-jf}
    \end{subfigure}

    \medskip

    \begin{subfigure}[t]{0.48\textwidth}
        \centering
        \includegraphics[width=\textwidth]
        {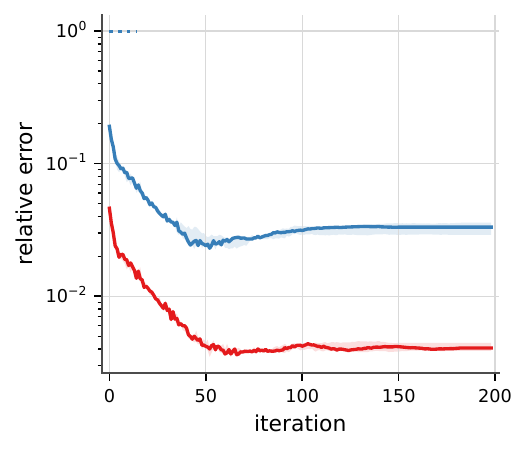}
        \captionsetup{justification=centering}
        \caption{LS: line source.%
        }
        \label{fig:exp-petrov-projection-ls}
    \end{subfigure}
    \hfill
    \begin{subfigure}[t]{0.48\textwidth}
        \centering
        \includegraphics[width=\textwidth]
        {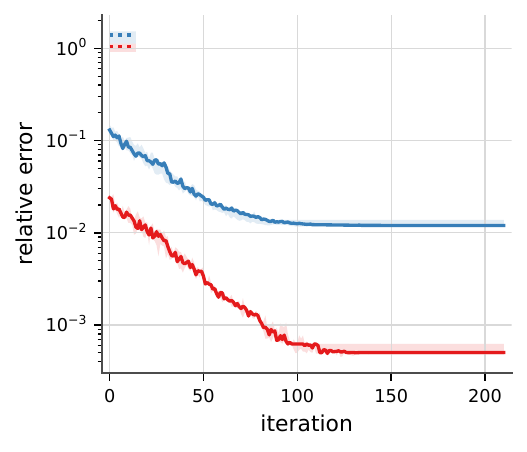}
        \caption{RC: reentrant corner.}
        \label{fig:exp-petrov-projection-rc}
    \end{subfigure}

    \medskip

    \includegraphics[width=0.72\textwidth]
    {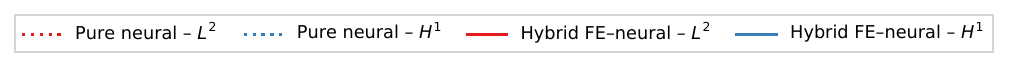}

    \caption{%
        \textbf{Effect of the hybrid finite element--neural construction for
        generic weak test functions.}
        Median relative $L^2$ and full relative $H^1$ errors over five seeds,
        with interquartile bands, for pure neural and hybrid finite
        element--neural models using the same compact random-hat test family.
        The hybrid construction is essentially neutral for MS, where the pure
        neural model already reaches very high accuracy, but substantially
        improves JF, LS, and RC.
        The pure neural runs for JF, LS, and RC terminate after $15$ iterations
        under the optimizer stopping criterion.
    }
    \label{fig:exp-petrov-projection}
\end{figure}

\begin{table}[t]
    \centering
    \small
    \caption{%
        Median final relative errors over five seeds.
        ``Best pure weak'' and ``best hybrid'' select the test family with the
        smallest relative $H^1$ error within the corresponding model class;
        the reported $L^2$ value is that of the same selected method.
        No strong-residual result is reported for LS because its line source
        is distributional. The JF finite element value is exceptional because
        the interface-fitted $P_4$ space contains the manufactured solution.
    }
    \label{tab:exp-general-summary}
    \begin{tabular}{llcccc}
        \toprule
        Method & Metric & MS & JF & LS & RC \\
        \midrule

        \multirow{2}{*}{GN Deep Ritz}
        & $L^2$
        & $4.48\times10^{1}$
        & $1.71\times10^{-3}$
        & $1.31\times10^{-2}$
        & $2.28\times10^{-1}$ \\
        & $H^1$
        & $5.25\times10^{2}$
        & $2.53\times10^{-2}$
        & $1.14\times10^{-1}$
        & $1.15$ \\

        \midrule

        \multirow{2}{*}{Strong PINN}
        & $L^2$
        & $9.98\times10^{-1}$
        & $5.97\times10^{-4}$
        & ---
        & $1.60\times10^{-3}$ \\
        & $H^1$
        & $9.98\times10^{-1}$
        & $1.60\times10^{-3}$
        & ---
        & $1.26\times10^{-2}$ \\

        \midrule

        \multirow{2}{*}{Best pure weak}
        & $L^2$
        & $2.13\times10^{-15}$
        & $1.45\times10^{-1}$
        & $1.16\times10^{-4}$
        & $1.05$ \\
        & $H^1$
        & $1.06\times10^{-13}$
        & $2.04$
        & $5.61\times10^{-3}$
        & $1.40$ \\

        \midrule

        \multirow{2}{*}{Best hybrid}
        & $L^2$
        & $1.49\times10^{-13}$
        & $4.26\times10^{-5}$
        & $6.72\times10^{-4}$
        & $5.04\times10^{-4}$ \\
        & $H^1$
        & $1.08\times10^{-12}$
        & $4.26\times10^{-4}$
        & $6.66\times10^{-3}$
        & $1.20\times10^{-2}$ \\

        \midrule

        \multirow{2}{*}{FE near hybrid budget}
        & $L^2$
        & $2.19\times10^{-7}$
        & $2.95\times10^{-14}$
        & $1.36\times10^{-3}$
        & $1.23\times10^{-3}$ \\
        & $H^1$
        & $7.31\times10^{-5}$
        & $5.44\times10^{-14}$
        & $1.48\times10^{-2}$
        & $1.60\times10^{-2}$ \\

        \bottomrule
    \end{tabular}
\end{table}

The results first show that no single test construction is uniformly
preferable. On JF, the best hybrid method uses tensor products of the
one-dimensional $H_0^1(0,1)$ sections and reaches relative errors
$4.26\times10^{-5}$ in $L^2$ and $4.26\times10^{-4}$ in $H^1$.
On RC, compact local hats are best, reaching
$5.04\times10^{-4}$ and $1.20\times10^{-2}$, respectively.
On LS, by contrast, the eigen-Green family gives the best pure neural result,
with errors $1.16\times10^{-4}$ in $L^2$ and
$5.61\times10^{-3}$ in $H^1$. Finally, on MS the pure
$H_0^1$ family reaches $2.13\times10^{-15}$ and
$1.06\times10^{-13}$. The test space is therefore a numerical design choice
rather than a fixed consequence of the operator.

The LS experiment illustrates a more structural advantage of the weak
formulation. Its line source defines a bounded functional on the variational
test space but cannot be represented by an ordinary pointwise residual.
Nevertheless, the weak solver reaches errors of order $10^{-4}$ in $L^2$ and
$10^{-3}$ in $H^1$. Together with the multiscale and reentrant-corner
benchmarks, this shows that the finite-dimensional Gauss--Newton construction
extends naturally to problems for which the strong residual is inconvenient
or not defined pointwise.

To isolate the effect of the finite element component in the hybrid
construction, we compare pure neural and hybrid finite element--neural
solvers using the same random-hat test family.
The corresponding final errors are reported in
\Cref{tab:exp-random-hats-pure-hybrid}.

\begin{table}[t]
    \centering
    \small
    \caption{%
        Median final relative errors over five seeds for pure neural and
        hybrid finite element--neural solvers using the same random-hat test family.
    }
    \label{tab:exp-random-hats-pure-hybrid}
    \begin{tabular}{llcccc}
        \toprule
        Method & Metric & MS & JF & LS & RC \\
        \midrule

        \multirow{2}{*}{Pure}
        & $L^2$
        & $2.19\times10^{-13}$
        & $2.05$
        & $9.93\times10^{-1}$
        & $1.052$ \\
        & $H^1$
        & $8.26\times10^{-12}$
        & $4.43$
        & $9.89\times10^{-1}$
        & $1.40$ \\

        \midrule

        \multirow{2}{*}{Hybrid}
        & $L^2$
        & $9.14\times10^{-13}$
        & $1.71\times10^{-4}$
        & $4.07\times10^{-3}$
        & $5.04\times10^{-4}$ \\
        & $H^1$
        & $2.16\times10^{-11}$
        & $1.34\times10^{-3}$
        & $3.32\times10^{-2}$
        & $1.20\times10^{-2}$ \\

        \bottomrule
    \end{tabular}
\end{table}

Hybridization is not a monotone improvement in final accuracy. On MS, the
pure neural solver already resolves the solution to very high accuracy, and
the hybrid construction provides essentially no further improvement at
convergence.
On JF, LS, and RC, however, it converts the same generic weak test
construction into an accurate solver. The relative $H^1$ errors decrease from
$4.43$ to $1.34\times10^{-3}$ on JF, from
$9.89\times10^{-1}$ to $3.32\times10^{-2}$ on LS, and from
$1.40$ to $1.20\times10^{-2}$ on RC, with analogous improvements in
$L^2$.

The standalone finite element (FE) references in
\Cref{tab:exp-general-summary} are selected from degree and mesh sweeps with
approximation dimensions between approximately $0.8$ and $1.2$ times the
corresponding hybrid budget.
At these comparable dimensions, the hybrid models substantially outperform
the selected finite element reference on MS, improve upon it on LS, and are
moderately more accurate on RC in both reported norms. 
JF is exceptional:
once the interface $x=1/2$ is fitted, the manufactured solution belongs to
the tested $P_4$ space and is therefore recovered to machine precision.
For RC, the selected finite element reference has $2409$ degrees of freedom,
compared with a total approximation dimension of $2008$ for the hybrid model.
These standalone references are independent of the coarse finite element
compensation spaces based on the meshes shown in
\Cref{fig:exp-petrov-compensation-meshes}.

To complement the comparison at fixed representation dimension, we also
report indicative wall-clock timings for the selected standalone finite
element references and for the corresponding pure neural and hybrid methods.
For the neural methods, we report both the median optimization time
$T_{\mathrm{optim}}$ and the median total runtime
$T_{\mathrm{total}}$ over the available seeds; for the finite element
references, we report the median warm solve time $T_{\mathrm{FE}}$.
These timings should be interpreted cautiously: they come from the present
implementations and hardware setup, and they measure computational cost
rather than approximation power. In particular, the finite element timings
concern standalone solves, whereas the neural timings correspond to iterative
optimization procedures.

\begin{table}[t]
    \centering
    \small
    \caption{%
        Indicative wall-clock timings for the pure neural, hybrid, and
        standalone finite element methods.
        For the neural methods, $T_{\mathrm{optim}}$ denotes the median
        optimization time and $T_{\mathrm{total}}$ the median total runtime.
        For the finite element references, $T_{\mathrm{FE}}$ denotes the
        median warm solve time.
    }
    \label{tab:exp-general-timings}
    \begin{tabular}{llccc}
        \toprule
        Benchmark & Method & $T_{\mathrm{optim}}$ (s)
        & $T_{\mathrm{total}}$ (s) & $T_{\mathrm{FE}}$ (s) \\
        \midrule

        \multirow{2}{*}{MS}
        & Pure neural & 299.5 & 304.2 & \multirow{2}{*}{0.040} \\
        & Hybrid      &  26.4 &  28.9 &                        \\
        \midrule

        \multirow{2}{*}{JF}
        & Pure neural &  75.7 &  80.0 & \multirow{2}{*}{0.309} \\
        & Hybrid      &  36.2 &  38.8 &                        \\
        \midrule

        \multirow{2}{*}{LS}
        & Pure neural & 134.5 & 137.3 & \multirow{2}{*}{0.345} \\
        & Hybrid      &  22.5 &  25.1 &                        \\
        \midrule

        \multirow{2}{*}{RC}
        & Pure neural$^\dagger$ & 9.7 & 12.3 & \multirow{2}{*}{2.172} \\
        & Hybrid                & 36.2 & 38.4 &                        \\

        \bottomrule
    \end{tabular}

    \vspace{0.3em}
    \begin{minipage}{0.94\textwidth}
        \footnotesize
        $^\dagger$The pure neural RC run stops after $15$ iterations without
        reaching an accurate solution; the corresponding runtime is therefore
        not representative of successful convergence.
    \end{minipage}
\end{table}

The timings in \Cref{tab:exp-general-timings} show that hybridization can
substantially reduce the computational cost of the neural solve. Relative to
the corresponding pure neural method, the total runtime is reduced by factors
of approximately $10.5$, $2.1$, and $5.5$ on MS, JF, and LS, respectively.
The apparently shorter pure neural runtime on RC should not be interpreted as
an efficiency advantage: the optimizer stalls after only $15$ iterations and
does not reach an accurate solution, whereas the hybrid method continues to
convergence.

Standalone finite element solves remain considerably faster in absolute
wall-clock time for the discretizations considered here. This comparison,
however, should be distinguished from representation efficiency. At comparable
approximation dimensions, the hybrid models are substantially more accurate
than the selected finite element reference on MS, improve upon it on LS, and
are slightly more accurate on RC, while JF is a favorable special case for
interface-fitted finite elements. Thus, the hybrid construction can both
reduce the optimization cost relative to a pure neural solver and provide a
more compact high-accuracy approximation than a finite element space of
comparable dimension.

The comparison is also relevant in light
of~\cite{grossmannCanPhysicsinformedNeural2024},
where standard PINN formulations were found to compare unfavorably with
finite element methods in both accuracy and computational cost. Subsequent
natural-gradient approaches
\cite{mullerAchievingHighAccuracy2023,SchF25,jniniGaussnewtonNaturalGradient2025,
schwencke2025amstramgram,jniniDualNaturalGradient2025,
jniniCurvatureAwareOptimizationHighAccuracy2026,
mckayNearoptimalSketchyNatural2025,
webbOptimisationFrameworkWellConditioned2026}
have already shown that the accuracy attainable by PINNs can depend strongly
on the optimization method. The present experiments extend this observation
to weak formulations: once the functional residual is discretized and
optimized within the Gauss--Newton framework developed above, weak neural
solvers can attain high accuracy even for problems for which pointwise
residual formulations are poorly suited or not defined. The hybrid finite
element--neural construction further improves the robustness of such
formulations for generic test functions.
A direct quantitative comparison with the benchmark study
of~\cite{grossmannCanPhysicsinformedNeural2024}, under matched problem,
implementation, and hardware conditions, is left for future work.

\subsubsection{Representative solution decompositions}

To complement the aggregate error measures, we examine representative
pointwise diagnostics for MS and RC. For each problem, we select the seed
whose final relative $H^1$ error is closest to the five-seed median of the
corresponding hybrid method; no best-seed selection is performed. For MS, we
also compare with an independently computed finite element approximation near
the hybrid budget.

\begin{figure}[t]
    \centering

    \begin{subfigure}[t]{0.48\textwidth}
        \centering
        \includegraphics[width=\textwidth]
        {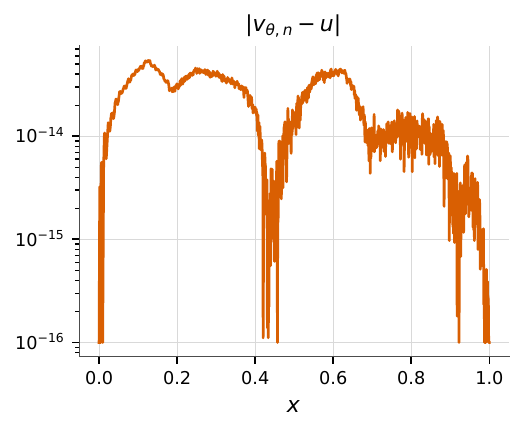}
        \caption{Hybrid error $|v_{\theta,n}-u|$.}
        \label{fig:exp-ms-hybrid-error}
    \end{subfigure}
    \hfill
    \begin{subfigure}[t]{0.48\textwidth}
        \centering
        \includegraphics[width=\textwidth]
        {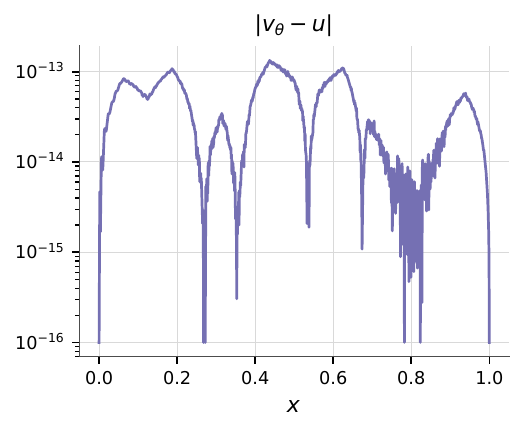}
        \caption{Neural-component error $|v_{\theta}-u|$.}
        \label{fig:exp-ms-neural-error}
    \end{subfigure}

    \medskip

    \begin{subfigure}[t]{0.48\textwidth}
        \centering
        \includegraphics[width=\textwidth]
        {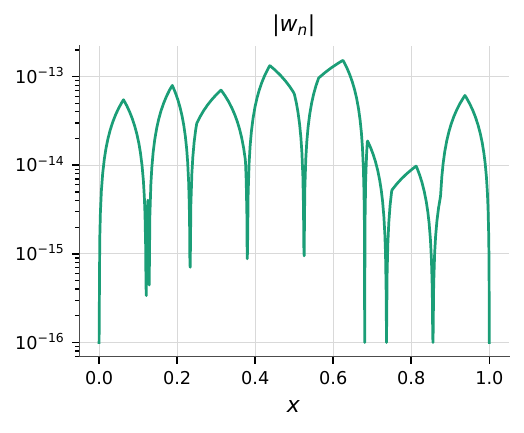}
        \caption{Magnitude of finite element compensator~$|w_n|$.}
        \label{fig:exp-ms-fe-compensator}
    \end{subfigure}
    \hfill
    \begin{subfigure}[t]{0.48\textwidth}
        \centering
        \includegraphics[width=\textwidth]
        {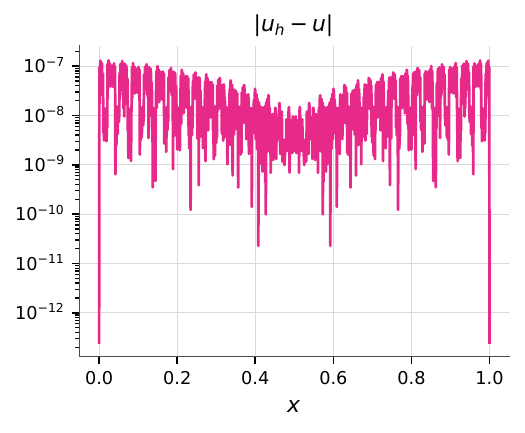}
        \caption{Finite element near hybrid budget error $|u_h-u|$.}
        \label{fig:exp-ms-fe-reference-error}
    \end{subfigure}

    \caption{%
        \textbf{MS: diagnostics for the hybrid decomposition.}
        The first three panels correspond to the representative hybrid run:
        the hybrid approximation error $|v_{\theta,n}-u|$, the neural
        component error $|v_{\theta}-u|$, and the magnitude $|w_n|$ of the
        finite element compensator. The fourth panel shows the pointwise absolute error
        of the independent finite element approximation near the hybrid budget. The
        selected finite element reference uses degree-$4$ elements on $96$ uniform
        elements, with $383$ degrees of freedom, compared with a total approximation
        dimension of $390$ for the hybrid model.
    }
    \label{fig:exp-ms-diagnostics}
\end{figure}

For MS, the representative diagnostics are shown in
\Cref{fig:exp-ms-diagnostics}. Panels
\subref{fig:exp-ms-hybrid-error}--\subref{fig:exp-ms-fe-compensator}
correspond to the representative hybrid run, while
\subref{fig:exp-ms-fe-reference-error} shows the independent finite element reference computed with a number of free parameters close to the hybrid model.

The diagnostics in \Cref{fig:exp-ms-diagnostics} show that the hybrid and
neural errors are essentially indistinguishable, while the finite element compensator
itself remains numerically negligible. Thus the hybrid approximation is
carried almost entirely by the neural component.

This behavior, while initially counter-intuitive, has a simple structural
explanation. Writing
\begin{equation*}
    v_{\theta,n}=v_\theta+w_n,
    \qquad
    w_n\in V_n,
\end{equation*}
exact recovery of the smooth solution would imply
\begin{equation*}
    w_n=u-v_\theta.
\end{equation*}
For the present experiment, both $u$ and $v_\theta$ are smooth, whereas
$w_n$ belongs to the continuous piecewise-affine finite element space $V_n$.
Hence exact recovery would require
\begin{equation*}
    w_n\in V_n\cap C^1([0,1]).
\end{equation*}
A globally $C^1$ piecewise-affine function is necessarily globally affine;
combined with the homogeneous Dirichlet conditions, this gives
\begin{equation*}
    V_n\cap C^1([0,1])=\{0\}.
\end{equation*}
Thus, in the exact-recovery limit, necessarily $w_n=0$ and
$v_\theta=u$. The numerically negligible finite element contribution observed in
\Cref{fig:exp-ms-diagnostics} is therefore consistent with the regularity
mismatch: introducing a nonzero piecewise-affine compensator would force the
neural component to approximate the less regular difference $u-w_n$, rather
than the smooth solution $u$ itself.

This argument, however, concerns the converged decomposition and not the
optimization path leading to it. During training, the finite element
component can make a
non-negligible contribution before becoming negligible near convergence,
consistently with the behavior observed in
\Cref{subsubsec:hybrid-results}. Thus, the nearly vanishing final compensator
should not be interpreted as indicating that hybridization does not play a role in 
the optimization. More detailed trajectory diagnostics are provided in
\weblink{https://nilo.schwencke.me/tutorials/beyond-pinns-companion/}{the companion blog}.

The finite element approximation with a budget of free parameters close to the hybrid model (shown in
\Cref{fig:exp-ms-diagnostics}) serves a different purpose. It is an
independent finite element reference at approximately the same total approximation
budget and should not be confused with the much coarser finite element compensation space with $15$ degrees of freedom used inside the hybrid model.

\begin{figure}[t]
    \centering
    \captionsetup[subfigure]{justification=centering}

    \begin{subfigure}[t]{0.31\textwidth}
        \centering
        \includegraphics[width=\textwidth]
        {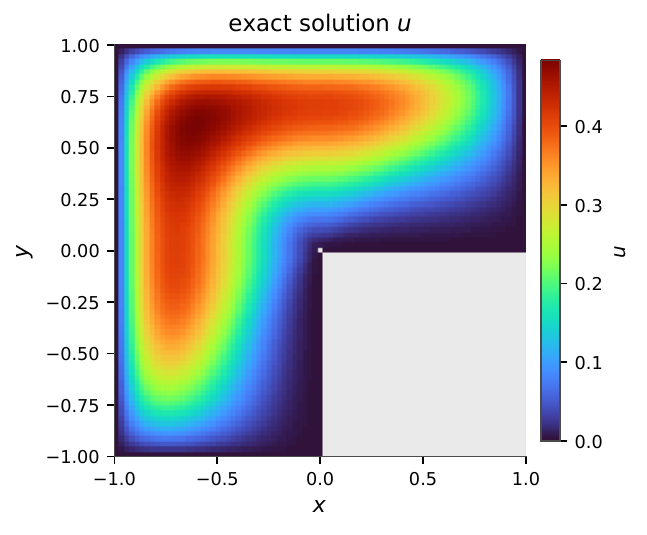}
        \caption{Exact solution $u$.}
        \label{fig:exp-rc-exact}
    \end{subfigure}
    \hfill
    \begin{subfigure}[t]{0.31\textwidth}
        \centering
        \includegraphics[width=\textwidth]
        {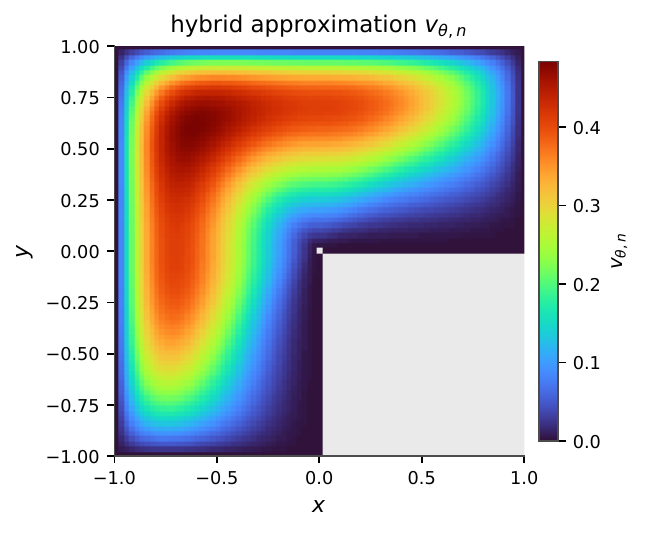}
        \caption{Hybrid approximation $v_{\theta,n}$.}
        \label{fig:exp-rc-hybrid}
    \end{subfigure}
    \hfill
    \begin{subfigure}[t]{0.31\textwidth}
        \centering
        \includegraphics[width=\textwidth]
        {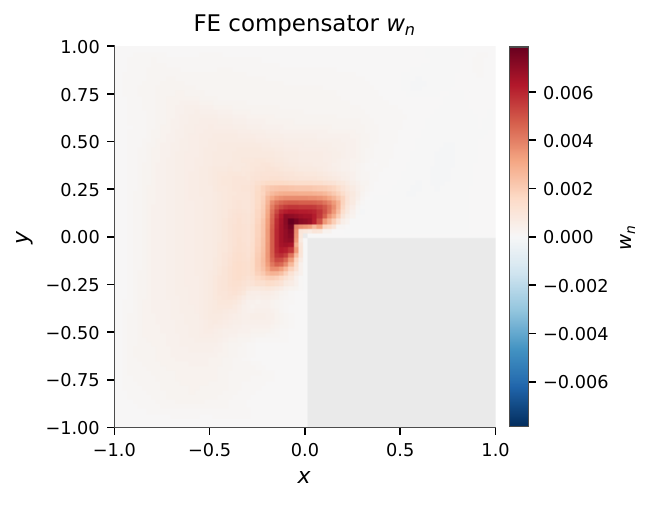}
        \caption{Finite element compensator $w_n$.}
        \label{fig:exp-rc-fe-compensator}
    \end{subfigure}

    \caption{%
        \textbf{RC: hybrid solution decomposition.}
        Exact solution $u$, hybrid approximation $v_{\theta,n}$, and finite
        element compensator $w_n$ for the representative hybrid run.
        The exact and hybrid solutions are displayed on the same color scale,
        whereas the finite element compensator uses its own scale to reveal its much
        smaller, localized contribution near the reentrant corner.
    }
    \label{fig:exp-rc-decomposition}
\end{figure}

For RC, the representative decomposition is shown in
\Cref{fig:exp-rc-decomposition}. In contrast to MS, the finite element compensator is
non-negligible at convergence and is strongly localized near the reentrant
corner.

\begin{figure}[t]
    \centering

    \begin{subfigure}[t]{\textwidth}
        \centering
        \includegraphics[width=\textwidth]
        {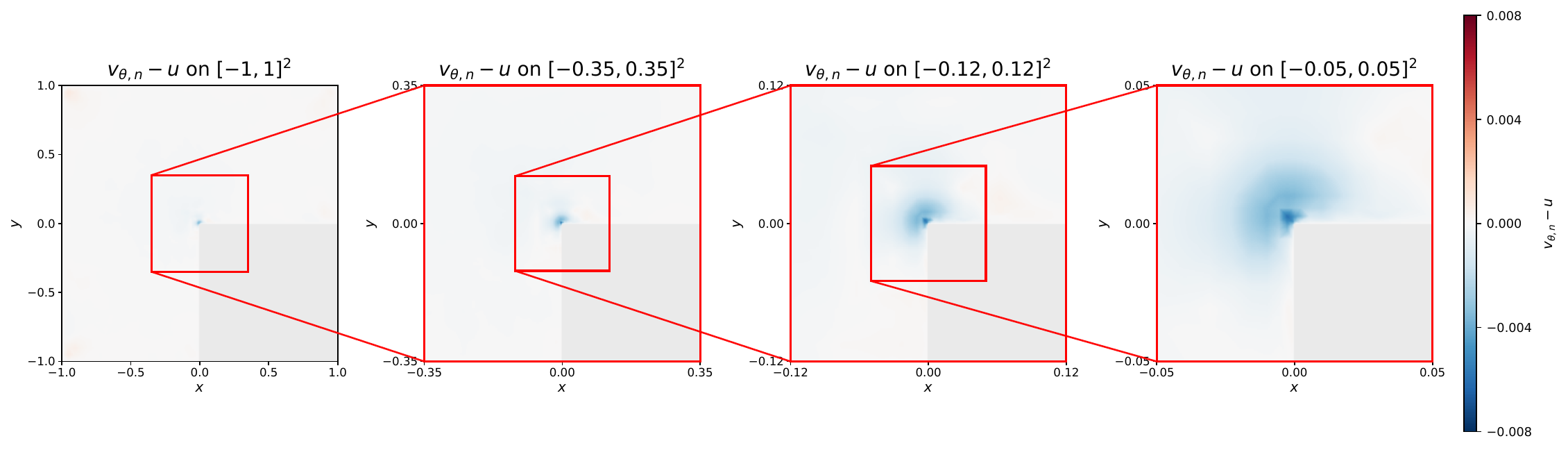}
        \caption{Signed hybrid difference $v_{\theta,n}-u$.}
        \label{fig:exp-rc-hybrid-difference}
    \end{subfigure}

    \medskip

    \begin{subfigure}[t]{\textwidth}
        \centering
        \includegraphics[width=\textwidth]
        {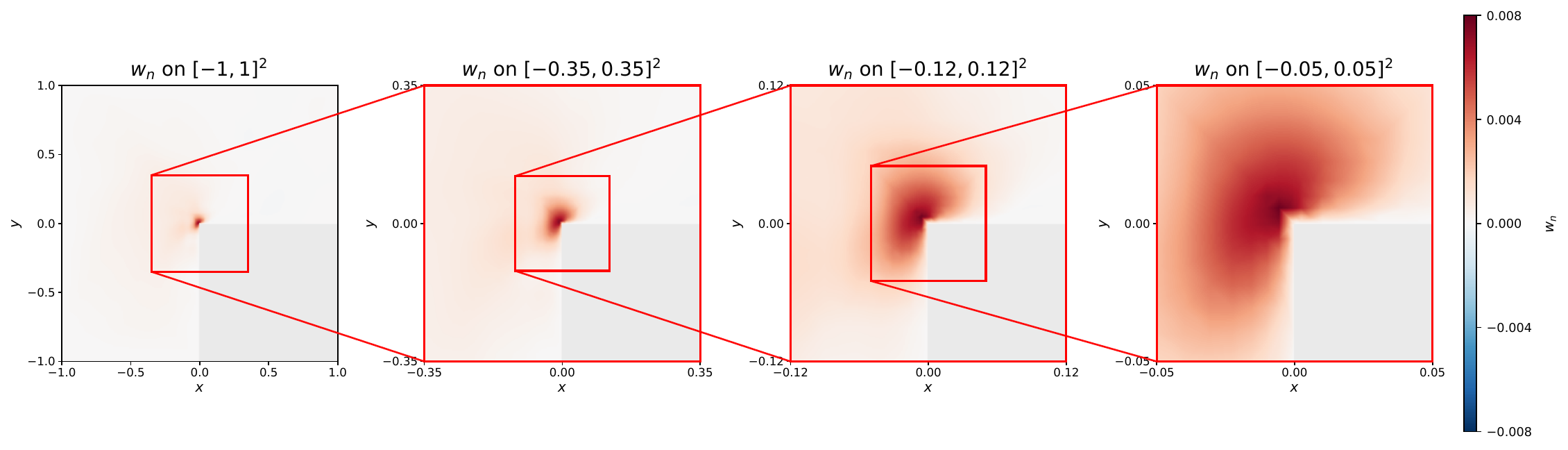}
        \caption{Signed finite element compensator $w_n$.}
        \label{fig:exp-rc-fe-compensator-corner}
    \end{subfigure}

    \caption{%
        \textbf{RC: localization near the reentrant corner.}
        Signed hybrid difference $v_{\theta,n}-u$ and finite element
        compensator $w_n$ for the representative hybrid run. From left to
        right, the panels show the full domain and the nested windows
        $[-0.35,0.35]^2$, $[-0.12,0.12]^2$, and
        $[-0.05,0.05]^2$, restricted to the L-shaped domain.
        Both quantities are displayed with the same symmetric color scale,
        allowing their magnitudes and spatial localization to be compared
        directly.
        The multilevel zoom layout and associated plotting implementation are adapted from~\cite{combetteNewInitialisationControl2026}.
    }
    \label{fig:exp-rc-corner}
\end{figure}

To examine this localization more closely,
\Cref{fig:exp-rc-corner} shows the signed hybrid difference
$v_{\theta,n}-u$ together with the signed finite element compensator $w_n$
over successively smaller neighborhoods of the corner.

Unlike the MS case, the finite element contribution does not vanish in the converged RC
approximation. Instead, \Cref{fig:exp-rc-fe-compensator-corner} shows a
localized, sign-changing correction concentrated around the reentrant corner,
where the exact solution has reduced regularity. The successive zooms show
that this structure persists down to the smallest displayed scale; the
$[-0.05,0.05]^2$ window still contains several layers of the graded finite element mesh.
The signed hybrid difference in
\Cref{fig:exp-rc-hybrid-difference} provides the corresponding local
approximation diagnostic on the same spatial and color scales.

Detailed test-family comparisons, finite element mesh and degree sweeps,
per-seed results, hybrid-update ablations, and further numerical
diagnostics are provided in
\weblink{https://nilo.schwencke.me/tutorials/beyond-pinns-companion/}{the companion blog}.

\bibliographystyle{abbrvnat}
\bibliography{bib}

\begin{thebibliography}{45}
\providecommand{\natexlab}[1]{#1}
\providecommand{\url}[1]{\texttt{#1}}
\expandafter\ifx\csname urlstyle\endcsname\relax
  \providecommand{\doi}[1]{doi: #1}\else
  \providecommand{\doi}{doi: \begingroup \urlstyle{rm}\Url}\fi

\bibitem[Barucq et~al.(2025)Barucq, Duprez, Faucher, Franck, Lecourtier,
  Lleras, {Michel-Dansac}, and Victorion]{BarDupFauFraLecLleMicVic2025}
H.~Barucq, M.~Duprez, F.~Faucher, E.~Franck, F.~Lecourtier, V.~Lleras,
  V.~{Michel-Dansac}, and N.~Victorion.
\newblock Enriching continuous {{Lagrange}} finite element approximation spaces
  using neural networks.
\newblock \emph{ArXiv Preprint}, 2502.04947, 2025.

\bibitem[Baydin et~al.(2018)Baydin, Pearlmutter, Radul, and
  Siskind]{baydinAutomaticDifferentiationMachine2018}
A.~G. Baydin, B.~A. Pearlmutter, A.~A. Radul, and J.~M. Siskind.
\newblock Automatic differentiation in machine learning: A survey.
\newblock \emph{Journal of Marchine Learning Research}, 18:\penalty0 1--43,
  2018.

\bibitem[Bon et~al.(2025)Bon, Caris, and Mula]{bonStableNonlinearDynamical2025}
D.~Bon, B.~Caris, and O.~Mula.
\newblock Stable {{Nonlinear Dynamical Approximation}} with {{Dynamical
  Sampling}}, May 2025.

\bibitem[Braess(2007)]{Bra07}
D.~Braess.
\newblock \emph{Finite Elements: Theory, Fast Solvers, and Applications in
  Solid Mechanics}.
\newblock Cambridge University Press, Cambridge, 2007.

\bibitem[Brenner and Scott(2008)]{BreS08}
S.~C. Brenner and L.~R. Scott.
\newblock \emph{The Mathematical Theory of Finite Element Methods}, volume~15
  of \emph{Texts in Applied Mathematics}.
\newblock Springer, New York, 3 edition, 2008.
\newblock \doi{10.1007/978-0-387-75934-0}.

\bibitem[Cho et~al.(2014)Cho, van Merrienboer, Bahdanau, and
  Bengio]{choPropertiesNeuralMachine2014}
K.~Cho, B.~van Merrienboer, D.~Bahdanau, and Y.~Bengio.
\newblock On the {{Properties}} of {{Neural Machine Translation}}:
  {{Encoder-Decoder Approaches}}.
\newblock \emph{ArXiv Preprint}, 1409.1259, 2014.
\newblock \doi{10.48550/arXiv.1409.1259}.

\bibitem[Ciarlet(1978)]{Cia78}
P.~G. Ciarlet.
\newblock \emph{The Finite Element Method for Elliptic Problems}, volume~4 of
  \emph{Studies in Mathematics and its Applications}.
\newblock North-Holland, Amsterdam, 1978.

\bibitem[Combette et~al.(2025)Combette, Venaille, and
  Pustelnik]{combetteNewInitialisationControl2026}
A.~Combette, A.~Venaille, and N.~Pustelnik.
\newblock A new initialisation to {{Control Gradients}} in {{Sinusoidal
  Neural}} network.
\newblock \emph{ArXiv Preprint}, 2512.06427, 2025.

\bibitem[Daw et~al.(2022)Daw, Bu, Wang, Perdikaris, and
  Karpatne]{dawMitigatingPropagationFailures2023}
A.~Daw, J.~Bu, S.~Wang, P.~Perdikaris, and A.~Karpatne.
\newblock Mitigating {{Propagation Failures}} in {{Physics-informed Neural
  Networks}} using {{Retain-Resample-Release}} ({{R3}}) {{Sampling}}, 2022.

\bibitem[Demkowicz and
  Gopalakrishnan(2011)]{demkowiczClassDiscontinuousPetrov2011}
L.~Demkowicz and J.~Gopalakrishnan.
\newblock A class of discontinuous {{Petrov}}--{{Galerkin}} methods. {{II}}.
  {{Optimal}} test functions.
\newblock \emph{Numerical Methods for Partial Differential Equations},
  27\penalty0 (1):\penalty0 70--105, 2011.

\bibitem[Dissanayake and
  Phan-Thien(1994)]{dissanayakeNeuralnetworkbasedApproximationsSolving1994}
M.~W. M.~G. Dissanayake and N.~Phan-Thien.
\newblock Neural-network-based approximations for solving partial differential
  equations.
\newblock \emph{Communications in Numerical Methods in Engineering},
  10\penalty0 (3):\penalty0 195--201, 1994.

\bibitem[Elman(1990)]{elmanFindingStructureTime1990}
J.~L. Elman.
\newblock Finding {{Structure}} in {{Time}}.
\newblock \emph{Cognitive Science}, 14\penalty0 (2):\penalty0 179--211, Mar.
  1990.

\bibitem[Franck et~al.(2026)Franck, {Michel-Dansac}, Navoret, and
  Vigon]{FraMicNavVig2025}
E.~Franck, V.~{Michel-Dansac}, L.~Navoret, and V.~Vigon.
\newblock Neural semi-{{Lagrangian}} method for high-dimensional
  advection-diffusion problems.
\newblock \emph{Computer Methods in Applied Mechanics and Engineering},
  448\penalty0 (B):\penalty0 118481, 2026.

\bibitem[Goodfellow et~al.(2016)Goodfellow, Bengio, and
  Courville]{goodfellowDeepLearning2016}
I.~Goodfellow, Y.~Bengio, and A.~Courville.
\newblock \emph{Deep Learning}.
\newblock MIT press, 2016.

\bibitem[Grossmann et~al.(2024)Grossmann, Komorowska, Latz, and
  Sch{\"o}nlieb]{grossmannCanPhysicsinformedNeural2024}
T.~G. Grossmann, U.~J. Komorowska, J.~Latz, and C.-B. Sch{\"o}nlieb.
\newblock Can physics-informed neural networks beat the finite element method?
\newblock \emph{IMA Journal of Applied Mathematics}, 89\penalty0 (1):\penalty0
  143--174, 2024.

\bibitem[Hao et~al.(2023)Hao, Hong, and Jin]{haoGaussNewtonVariational2023}
W.~Hao, Q.~Hong, and X.~Jin.
\newblock Gauss {N}ewton method for solving variational problems of {{PDEs}}
  with neural network discretizaitons.
\newblock \emph{ArXiv Preprint}, 2306.08727, 2023.

\bibitem[Hochreiter and Schmidhuber(1997)]{hochreiterLongShorttermMemory1997}
S.~Hochreiter and J.~Schmidhuber.
\newblock Long short-term memory.
\newblock \emph{Neural computation}, 9\penalty0 (8):\penalty0 1735--1780, 1997.

\bibitem[Jnini and Vella(2025)]{jniniDualNaturalGradient2025}
A.~Jnini and F.~Vella.
\newblock Dual {{Natural Gradient Descent}} for {{Scalable Training}} of
  {{Physics-Informed Neural Networks}}.
\newblock \emph{ArXiv Preprint}, 2505.21404, 2025.

\bibitem[Jnini et~al.(2025)Jnini, Vella, and
  Zeinhofer]{jniniGaussnewtonNaturalGradient2025}
A.~Jnini, F.~Vella, and M.~Zeinhofer.
\newblock Gauss-{N}ewton natural gradient descent for physics-informed
  computational fluid dynamics.
\newblock \emph{Computers \& Fluids}, 307:\penalty0 106955, 2025.

\bibitem[Jnini et~al.(2026)Jnini, Kiyani, Shukla, Urban, Daryakenari, Muller,
  Zeinhofer, and Karniadakis]{jniniCurvatureAwareOptimizationHighAccuracy2026}
A.~Jnini, E.~Kiyani, K.~Shukla, J.~F. Urban, N.~A. Daryakenari, J.~Muller,
  M.~Zeinhofer, and G.~E. Karniadakis.
\newblock Curvature-{{Aware Optimization}} for {{High-Accuracy Physics-Informed
  Neural Networks}}.
\newblock \emph{ArXiv Preprint}, 2604.05230, 2026.

\bibitem[Kharazmi et~al.(2019)Kharazmi, Zhang, and
  Karniadakis]{kharazmiVariationalPhysicsInformedNeural2019}
E.~Kharazmi, Z.~Zhang, and G.~E. Karniadakis.
\newblock Variational {{Physics-Informed Neural Networks For Solving Partial
  Differential Equations}}.
\newblock \emph{ArXiv Preprint}, 1912.00873, 2019.

\bibitem[Kharazmi et~al.(2021)Kharazmi, Zhang, and
  Karniadakis]{kharazmiHpVPINNsVariationalPhysicsinformed2021}
E.~Kharazmi, Z.~Zhang, and G.~E. Karniadakis.
\newblock Hp-{{VPINNs}}: {{Variational}} physics-informed neural networks with
  domain decomposition.
\newblock \emph{Computer Methods in Applied Mechanics and Engineering},
  374:\penalty0 113547, 2021.

\bibitem[Lagaris et~al.(1998)Lagaris, Likas, and
  Fotiadis]{lagarisArtificialNeuralNetworks1998}
I.~E. Lagaris, A.~Likas, and D.~I. Fotiadis.
\newblock Artificial neural networks for solving ordinary and partial
  differential equations.
\newblock \emph{IEEE transactions on neural networks}, 9\penalty0 (5):\penalty0
  987--1000, 1998.

\bibitem[Lau et~al.(2024)Lau, Hemachandra, Ng, and
  Low]{lauPINNACLEPINNAdaptive2024a}
G.~K.~R. Lau, A.~Hemachandra, S.-K. Ng, and B.~K.~H. Low.
\newblock {{PINNACLE}}: {{PINN}} adaptive {{ColLocation}} and experimental
  points selection.
\newblock In \emph{The Twelfth International Conference on Learning
  Representations}, 2024.

\bibitem[LeCun et~al.(1998)LeCun, Bottou, Bengio, and
  Haffner]{lecunGradientbasedLearningApplied1998}
Y.~LeCun, L.~Bottou, Y.~Bengio, and P.~Haffner.
\newblock Gradient-based learning applied to document recognition.
\newblock \emph{Proceedings of the IEEE}, 86\penalty0 (11):\penalty0
  2278--2324, 1998.

\bibitem[Linnainmaa(1976)]{linnainmaaTaylorExpansionAccumulated1976}
S.~Linnainmaa.
\newblock Taylor expansion of the accumulated rounding error.
\newblock \emph{BIT Numerical Mathematics}, 16\penalty0 (2):\penalty0 146--160,
  1976.

\bibitem[M{\aa}lqvist and Peterseim(2021)]{MalP21}
A.~M{\aa}lqvist and D.~Peterseim.
\newblock \emph{Numerical homogenization by localized orthogonal
  decomposition}, volume~5 of \emph{SIAM Spotlights}.
\newblock Society for Industrial and Applied Mathematics (SIAM), Philadelphia,
  PA, 2021.

\bibitem[Mao and Meng(2023)]{maoPhysicsinformedNeuralNetworks2023}
Z.~Mao and X.~Meng.
\newblock Physics-informed neural networks with residual/gradient-based
  adaptive sampling methods for solving partial differential equations with
  sharp solutions.
\newblock \emph{Applied Mathematics and Mechanics}, 44\penalty0 (7):\penalty0
  1069--1084, 2023.

\bibitem[Margenberg et~al.(2024)Margenberg, Jendersie, Lessig, and
  Richter]{margenbergDNNMGHybridNeural2024}
N.~Margenberg, R.~Jendersie, C.~Lessig, and T.~Richter.
\newblock {{DNN-MG}}: {{A}} hybrid neural network/finite element method with
  applications to {{3D}} simulations of the {{Navier}}--{{Stokes}} equations.
\newblock \emph{Computer Methods in Applied Mechanics and Engineering},
  420:\penalty0 116692, 2024.

\bibitem[McKay et~al.(2025)McKay, Kaur, Greif, and
  Wetton]{mckayNearoptimalSketchyNatural2025}
M.~B. McKay, A.~Kaur, C.~Greif, and B.~Wetton.
\newblock Near-optimal {{Sketchy Natural Gradients}} for {{Physics-Informed
  Neural Networks}}, 2025.
\newblock URL \url{https://openreview.net/forum?id=bKsZomnmqn}.

\bibitem[McKay et~al.(2026)McKay, Lawrence, Wetton, and
  Gopaluni]{mckayErrorWhiteningWhy2026}
M.~B. McKay, N.~P. Lawrence, B.~Wetton, and R.~B. Gopaluni.
\newblock Error whitening: {{Why Gauss-Newton}} outperforms {{Newton}}.
\newblock \emph{ArXiv Preprint}, 2605.11316, 2026.

\bibitem[Müller and Zeinhofer(2023)]{mullerAchievingHighAccuracy2023}
J.~Müller and M.~Zeinhofer.
\newblock Achieving high accuracy with {{PINNs}} via energy natural gradient
  descent.
\newblock In \emph{International {{Conference}} on {{Machine Learning}}}, pages
  25471--25485. PMLR, 2023.

\bibitem[Nabian et~al.(2021)Nabian, Gladstone, and
  Meidani]{nabianEfficientTrainingPhysicsinformed2021}
M.~A. Nabian, R.~J. Gladstone, and H.~Meidani.
\newblock Efficient training of physics-informed neural networks via importance
  sampling.
\newblock \emph{Computer-Aided Civil and Infrastructure Engineering},
  36\penalty0 (8):\penalty0 962--977, 2021.

\bibitem[Nguyen et~al.(2023)Nguyen, Dairay, Meunier, Millet, and
  Mougeot]{nguyenFixedBudgetOnlineAdaptive2023}
T.~N.~K. Nguyen, T.~Dairay, R.~Meunier, C.~Millet, and M.~Mougeot.
\newblock Fixed-{{Budget Online Adaptive Learning}} for~{{Physics-Informed
  Neural Networks}}. {{Towards Parameterized Problem Inference}}.
\newblock In J.~Miky{\v s}ka, C.~{de Mulatier}, M.~Paszynski, V.~V.
  Krzhizhanovskaya, J.~J. Dongarra, and P.~M. Sloot, editors,
  \emph{Computational {{Science}} -- {{ICCS}} 2023}, pages 453--468, Cham,
  2023. Springer Nature Switzerland.

\bibitem[Nouy and Somacal(2026)]{nouyNaturalGradientDescent2026}
A.~Nouy and A.~Somacal.
\newblock Natural gradient descent with momentum.
\newblock \emph{ArXiv Preprint}, 2604.15554, 2026.

\bibitem[Paulsen and
  Raghupathi(2016)]{paulsenIntroductionTheoryReproducing2016}
V.~I. Paulsen and M.~Raghupathi.
\newblock \emph{An {{Introduction}} to the {{Theory}} of {{Reproducing Kernel
  Hilbert Spaces}}}, volume 152.
\newblock Cambridge University Press, 2016.

\bibitem[Raissi et~al.(2019)Raissi, Perdikaris, and
  Karniadakis]{raissiPhysicsinformedNeuralNetworks2019}
M.~Raissi, P.~Perdikaris, and G.~Karniadakis.
\newblock Physics-informed neural networks: {{A}} deep learning framework for
  solving forward and inverse problems involving nonlinear partial differential
  equations.
\newblock \emph{Journal of Computational Physics}, 378:\penalty0 686--707,
  2019.

\bibitem[Rojas et~al.(2024)Rojas, Maczuga, {Mu{\~n}oz-Matute}, Pardo, and
  Paszy{\'n}ski]{rojasRobustVariationalPhysicsInformed2024}
S.~Rojas, P.~Maczuga, J.~{Mu{\~n}oz-Matute}, D.~Pardo, and M.~Paszy{\'n}ski.
\newblock Robust {{Variational Physics-Informed Neural Networks}}.
\newblock \emph{Computer Methods in Applied Mechanics and Engineering},
  425:\penalty0 116904, 2024.

\bibitem[Schwencke and Furtlehner(2025)]{SchF25}
N.~Schwencke and C.~Furtlehner.
\newblock {{ANaGRAM}}: A natural gradient relative to adapted model for
  efficient {{PINNs}} learning.
\newblock In \emph{The Thirteenth International Conference on Learning
  Representations}, 2025.

\bibitem[Schwencke et~al.(2025)Schwencke, Rousselot, Shilova, and
  Furtlehner]{schwencke2025amstramgram}
N.~Schwencke, C.~Rousselot, A.~Shilova, and C.~Furtlehner.
\newblock {{AMStramGRAM}}: {{Adaptive}} multi-cutoff strategy modification for
  {{ANaGRAM}}.
\newblock \emph{ArXiv Preprint}, 2510.15998, 2025.

\bibitem[Shang et~al.(2022)Shang, Wang, and
  Sun]{shangDeepPetrovGalerkinMethod2022}
Y.~Shang, F.~Wang, and J.~Sun.
\newblock Deep {{Petrov-Galerkin Method}} for {{Solving Partial Differential
  Equations}}.
\newblock \emph{ArXiv Preprint}, 2201.12995, 2022.

\bibitem[Vaswani et~al.(2017)Vaswani, Shazeer, Parmar, Uszkoreit, Jones, Gomez,
  Kaiser, and Polosukhin]{vaswaniAttentionAllYou2017}
A.~Vaswani, N.~Shazeer, N.~Parmar, J.~Uszkoreit, L.~Jones, A.~N. Gomez,
  L.~Kaiser, and I.~Polosukhin.
\newblock Attention {{Is All You Need}}.
\newblock \emph{ArXiv Preprint}, 1706.03762, 2017.

\bibitem[Webb et~al.(2026)Webb, Jerad, and
  Cartis]{webbOptimisationFrameworkWellConditioned2026}
J.~Webb, S.~Jerad, and C.~Cartis.
\newblock An {{Optimisation Framework}} for the {{Well-Conditioned Training}}
  of {{Physics-Informed Neural Networks}}.
\newblock \emph{ArXiv Preprint}, 2607.02194, 2026.

\bibitem[Wu et~al.(2023)Wu, Zhu, Tan, Kartha, and
  Lu]{wuComprehensiveStudyNonadaptive2023}
C.~Wu, M.~Zhu, Q.~Tan, Y.~Kartha, and L.~Lu.
\newblock A comprehensive study of non-adaptive and residual-based adaptive
  sampling for physics-informed neural networks.
\newblock \emph{Computer Methods in Applied Mechanics and Engineering},
  403:\penalty0 115671, 2023.

\bibitem[Zang et~al.(2020)Zang, Bao, Ye, and
  Zhou]{zangWeakAdversarialNetworks2020}
Y.~Zang, G.~Bao, X.~Ye, and H.~Zhou.
\newblock Weak adversarial networks for high-dimensional partial differential
  equations.
\newblock \emph{Journal of Computational Physics}, 411:\penalty0 109409, 2020.

\end{thebibliography}

\appendix

\section{Classical methods as instances of the Petrov--Galerkin framework}
\label{app:gn-approximation-and-test-examples}

The framework of \Cref{subsec:discretized-functional-gn} separates two
ingredients: the local approximation space given by the tangent space of the
parametric model and the measurements, or equivalently test functions, used
to discretize the linearized residual equation.
We briefly record how several familiar
constructions arise from particular choices of these ingredients. Throughout,
we use the notation introduced in \Cref{subsec:discretized-functional-gn}.
For more detailed developments and complementary explanations, particularly
from a machine-learning perspective, we refer to \weblink{https://nilo.schwencke.me/tutorials/beyond-pinns-companion/}{the companion blog}.

\subsection{Natural gradient.}
Let $\pmFun:\R^p\to\calH$ be a parametric model with residual
$\resFun_\paramPM=\pmFun_\paramPM-y$. Choosing the tangent directions
$\partial_i\pmFun_\paramPM$ themselves as test functions in
\eqref{eq:functional-petrov-galerkin} gives
\begin{equation}
    \sum_{j=1}^p
    \braket{
        \partial_j\pmFun_\paramPM
    }{
        \partial_i\pmFun_\paramPM
    }_{\calH}
    \paramUpdate_j
    =
    \braket{
        \resFun_\paramPM
    }{
        \partial_i\pmFun_\paramPM
    }_{\calH},
    \qquad
    1\leq i\leq p.
    \label{eq:app-natural-gradient}
\end{equation}
Writing
$
G_{\paramPM,ij}
=
\braket{\partial_i\pmFun_\paramPM}
       {\partial_j\pmFun_\paramPM}_{\calH},
$
this is
$
G_\paramPM\paramUpdate=\nabla_\paramPM\ell(\paramPM)
$
for
$
\ell(\paramPM)=\frac12\|\resFun_\paramPM\|_{\calH}^2.
$
Hence the corresponding correction is the natural-gradient direction
\cite{mullerAchievingHighAccuracy2023,
jniniGaussnewtonNaturalGradient2025,
mckayErrorWhiteningWhy2026}, and its functional action satisfies
\begin{equation}
    \dd\pmFun_\paramPM[\paramUpdate]
    =
    \Pi_{T_\paramPM\calM_\pmFun}\resFun_\paramPM.
\end{equation}
Thus, natural gradient is the Galerkin choice in which the approximation and
test spaces both coincide with the tangent space.

\subsection{Kernel and empirical natural-gradient methods.}
Assume that point evaluation is well defined and continuous on a Hilbert
space $\calH$. Its Riesz representer $k(x,\bullet)\in\calH$ satisfies
\begin{equation}
    w(x)
    =
    \braket{
        w
    }{
        k(x,\bullet)
    }_{\calH}.
\end{equation}
For points $X=\{x_i\}_{i=1}^N$, choosing
$
\widehat{\calH}_X
=
\spn\{k(x_i,\bullet)\}_{i=1}^N
$
as both approximation and test space therefore yields the orthogonal
projection of $y$ onto $\widehat{\calH}_X$, equivalently the standard kernel
interpolant associated with the points $X$
\cite{paulsenIntroductionTheoryReproducing2016}.

The empirical natural gradient of \cite{SchF25} is the local version of this
construction. Point evaluation on the tangent space
$T_\paramPM\calM_\pmFun$ admits representers
$k_\paramPM(x_i,\bullet)$, and the corresponding empirical tangent space is
\begin{equation}
    \widehat T_{\paramPM,X}\calM_\pmFun
    =
    \spn\{
        k_\paramPM(x_i,\bullet)
        \mid 1\leq i\leq N
    \}
    \subset
    T_\paramPM\calM_\pmFun.
\end{equation}
The empirical natural-gradient correction is consequently characterized by
\begin{equation}
    \widehat{\pmFun}_{\paramPM,X}(\paramLinear)
    =
    \Pi_{\widehat T_{\paramPM,X}\calM_\pmFun}
    \resFun_\paramPM.
\end{equation}
It thus replaces the full tangent-space projection of natural gradient by an
orthogonal projection onto the kernel space generated by the chosen
measurements.

\subsection{Pointwise Gauss--Newton and PINN collocation.}
For a compound parametric model
$\pmComp:\R^p\to\calY$, let point evaluation at
$X=\{x_i\}_{i=1}^N$ be well defined on the corresponding extended tangent
space. Choosing
\begin{equation}
    \linform_i^\paramPM(w)=w(x_i)
\end{equation}
in \eqref{eq:discrete-functional-system} gives
\begin{equation}
    \resDis_{\paramPM,i}^{X}
    =
    \resFun_\paramPM(x_i),
    \qquad
    \jacDis_{\paramPM,ij}^{X}
    =
    \partial_j\pmComp_\paramPM(x_i).
    \label{eq:app-pointwise-gn}
\end{equation}
The resulting least-squares problem is precisely the standard pointwise
Gauss--Newton discretization. Taking $\pmComp$ to be the PDE compound model
therefore recovers the usual collocation-based Gauss--Newton formulation for
PINNs. Since point evaluation admits a Riesz representer on the
finite-dimensional extended tangent space, this construction is itself a
particular Petrov--Galerkin discretization of the form
\eqref{eq:functional-petrov-galerkin}.

Observational data are incorporated in the same way by augmenting the PDE
compound model with the identity, $(D,B,\Id[\calH])\circ\pmFun$, and applying
point-evaluation measurements to this additional component, thereby recovering
the usual supervised data term in PINNs.

\subsection{Sketching and regularization.}
The same viewpoint also gives concise interpretations of common algebraic
modifications of Gauss--Newton systems. If
$
P\in\R^{p\times s_p}
$
restricts parameter corrections and
$
S\in\R^{s_t\times N}
$
combines the measurements, the discretized system becomes
\begin{equation}
    S\jacDis_\paramPM^{\Lambda_\paramPM}P\,\eta
    =
    S\resDis_\paramPM^{\Lambda_\paramPM}.
\end{equation}
Hence parameter-side sketching restricts the local approximation space,
whereas residual-side sketching replaces the test family by linear
combinations of its elements
\cite{mckayNearoptimalSketchyNatural2025,
webbOptimisationFrameworkWellConditioned2026}.

Likewise, ridge-regularized Gauss--Newton follows by augmenting the functional
model with the scaled identity,
$
\pmFun^\alpha_\paramPM=(\pmFun_\paramPM,\alpha\paramPM).
$
After recentering the augmented target at the current iterate, the local
problem becomes
\begin{equation}
    \min_{\paramUpdate\in\R^p}
    \frac12
    \|
        \resFun_\paramPM-\dd\pmFun_\paramPM[\paramUpdate]
    \|_{\calH}^2
    +
    \frac{\alpha^2}{2}
    \|\paramUpdate\|_{\R^p}^2,
\end{equation}
which yields the usual ridge-regularized Gauss--Newton correction~\cite{schwencke2025amstramgram}.

\end{document}